\documentclass[a4paper,12p]{article}
\usepackage{amsthm,amsmath,amssymb,mathrsfs,url,amsfonts} 
\theoremstyle{definition}

\begin{document}

\title{Dynamical Bulk Scaling limit of Gaussian Unitary Ensembles and 
Stochastic Differential Equation gaps}
\author{Yosuke Kawamoto ; Hirofumi Osada \\
Appeared in J,Theoretical Probability: https://doi.org/10.1007/s10959-018-0816-2}    
\maketitle

\pagestyle{myheadings}%\pagestyle{plain}
\markboth{Dynamical Bulk Scaling Limit of GUE and SDE gaps}
{Yosuke Kawamoto, Hirofumi Osada}
%\markboth{Dynamical Bulk Scaling Limit of GUE and SDE gaps}

%\author{Yosuke Kawamoto,\quad Hirofumi Osada}

%\subjclass[2010]{Primary 60B20; Secondary 60H10}
\vskip 0.5cm

\newcommand\DN{\newcommand}\newcommand\DR{\renewcommand}

\newcommand{\ep}{\varepsilon}

 \theoremstyle{definition}
 \newtheorem{theorem}{Theorem}[section]
 \newtheorem*{theorem*}{Theorem}
 \newtheorem{lemma}[theorem]{Lemma}
 \newtheorem{proposition}[theorem]{Proposition}
 \newtheorem{corollary}[theorem]{Corollary}
 \newtheorem{definition}{Definition}[section]

\DN\lref[1]{Lemma~\ref{#1}}
\DN\tref[1]{Theorem~\ref{#1}}
\DN\pref[1]{Proposition~\ref{#1}}
\DN\sref[1]{Section~\ref{#1}}
\DN\dref[1]{Definition~\ref{#1}}
\DN\rref[1]{Remark~\ref{#1}} 
\DN\corref[1]{Corollary~\ref{#1}}
\DN\eref[1]{Example~\ref{#1}}

\numberwithin{equation}{section}
\numberwithin{theorem}{section}
% \numberwithin{lemma}{section}
% \numberwithin{definition}{section}
% \numberwithin{corollary}{section}
% \numberwithin{example}{section}

\newcounter{Const} \setcounter{Const}{0}
\DN\Ct[1]{\refstepcounter{Const}c_{\theConst}\label{#1}}\DN\cref[1]{c_{\ref{#1}}}

%極限
\DN{\limi}[1]{\lim_{#1\to\infty}} 	
\DN{\limz}[1]{\lim_{#1\to0}}
\DN{\limdz}[1]{\lim_{#1\downarrow 0}}
\DN{\limsupi}[1]{\limsup_{#1\to\infty}}
\DN{\limsupz}[1]{\limsup_{#1\to0}}
\DN{\limsupdz}[1]{\limsup_{#1\downarrow0}}
\DN{\liminfi}[1]{\liminf_{#1\to\infty}}
\DN{\liminfz}[1]{\liminf_{#1\to0}}
\DN{\liminfdz}[1]{\liminf_{#1\downarrow0}}
\DN{\supnor}[1]{\| #1\|_{\infty}}
\DN{\sumii}[1]{\sum_{#1=1}^{\infty}}
\DN{\sumi}[1]{\sum_{#1=0}^{\infty}}
\DN\tolaw{\stackrel{d}{\longrightarrow }}

%関数空間
\DN\Lm{L^{2}(\mu )}
\DN\Lploc{L_{\mathrm{loc}}^{p}}
\DN\Lqloc{L_{\mathrm{loc}}^{q}}
\DN\Loneloc{L_{\mathrm{loc}}^{1}}
\DN\Ltwoloc{L_{\mathrm{loc}}^{2}}
%\DN\Ci{C^{\infty}}
%\DN\Czi{C_{0}^{\infty}}

\DN\ord{{\mathcal o}}
\DN\Ord{\mathcal{O}}
\DN\PD[2]{\frac{\partial #1 }{\partial #2}}
\DN\half{\frac{1}{2}}
\DN\map[3]{#1:#2 \to #3}

\DN\R{\mathbb{R}}
\DN\Rd{\mathbb{R}^d}
\DN\N{\mathbb{N}}
\DN\Q{\mathbb{Q}}
\DN{\Z}{\mathbb{Z}}
\DN\C{\mathbb{C}}

\DN\bs{\bigskip}
\DN\ms{\medskip}
\DN\ts{\times}
\DN\bra{\langle}
\DN\ket{\rangle}

\DN\SSup[1]{\Big\| #1 \Big\|_{R}}
\DN\Sup[1]{\Big\| #1 \Big\|_{R}}

\DN\PROOF{\noindent{\it Proof. }}
\DN\PROOFEND{\qed}
%\DN\PROOF{\begin{proof}}
%\DN\PROOFEND{\qed \end{proof}}

%\input{sg-def}

\DN\aD{-\frac{2}{3}}
\DN\aU{-\half }

\DN\ff{\mathsf{g}}
\DN\nablax{\nabla_x}
\DN\p{\mathsf{p}}
\DN\phat{\hat{p}}

\DN\mul{\mu \circ \lab ^{-1}}
\DN\Sr{S_r}
\DN\SR{S_R}
\DN\sS{\R }
\DN\SSS{\mathsf{S}}
\DN\sss{\mathsf{s}}
\DN\SSSsi{\SSS _{\mathrm{s.i.}}}
\DN\SoneSSS{\sS \times \SSS }
\DN\SSSksingle{\SSSsi ^{k}}
\DN\Srm{\SSS _{r}^{m}}
\DN\SSSS{\mathbb{S}}

\DN\BN{\mathsf{B}^{N }}

\DN\U{\mathsf{U}}
\DN\UN{\U ^{N }}
\DN\UNN{N \UN } 
\DN\UNone{\UN _1}
\DN\UNNone{\UNN _1}
\DN\UNtwo{\U _2^N }

\DN\KN{\mathsf{K}^{\nN }}
\DN\KtN{\mathsf{K}_{\theta }^{\nN }}
\DN\nN{N }

\DN\xN{x_N}
\DN\yN{y_N}
\DN\zN{z_N}
\DN\tN{t_N}

\DN\xxN{\mathsf{x}_{\theta }^N}
\DN\xxNN{\mathsf{x}_{-\theta }^N}

\DN\LNx{\mathsf{L}^{\nN }(x,x)}
\DN\LNy{\mathsf{L}^{\nN }(y,y)}
\DN\LNyz{\mathsf{L}^{\nN }(y,z)}
\DN\LNzx{\mathsf{L}^{\nN }(z,x)}
\DN\LNzy{\mathsf{L}^{\nN }(z,y)}
\DN\LNNyz{\mathsf{L}^{\nN -1}(\zN ,\zN )}
\DN\LNyx{\mathsf{L}^{\nN }(\xxN ,y)}
\DN\LNxz{\mathsf{L}^{\nN }(\xxN ,z)}
\DN\LNxy{\mathsf{L}^{\nN }(x ,y)}
\DN\LNxyN{\mathsf{L}^{\nN }(\xxN ,y)}
\DN\LNz{\mathsf{L}^{\nN }(z,z)}
\DN\LNzN{\mathsf{L}^{\nN }(\zN ,\zN )}
\DN\LNxN{\mathsf{L}^{\nN }(\xxN ,\xxN )}
\DN\LNyxN{\mathsf{L}^{\nN }(\xxN ,y )}
\DN\LNzxN{\mathsf{L}^{\nN }(\xxN ,z )}
\DN\LNxNy{\mathsf{L}^{\nN }(\xxN ,y )}

\DN\ellxy{\xxN - y }
\DN\ellyz{y -z }
\DN\ellxz{ \xxN - z }

\DN\ellxyP{(\ellxy )}
\DN\ellyzP{(\ellyz )}
\DN\ellxzP{(\ellxz )}

\DN\ellxyA{ |\xxN - y |}
\DN\ellyzA{ |y - z |}
\DN\ellxzA{ |\xxN - z |}

\DN\ellXYZ{\ellxyA \ellxzA }
\DN\ellxyz{\ellxyP \ellxzP }

%一般記号
\DN\q{1}
\DN\qN{\frac{\q }{ N }}
\DN\oneN{1/N}

\DN\ww{\rr } %[]
\DN\www{\rr (x)}
\DN\rr{{w}}
\DN\rrN{\rr ^{N}}
\DN\rrNs{w_{r}}
\DN\rrbar{\bar{\rr }}

\DN\As[1]{\textbf{(#1)}}

%RPF
\DN\muone{\mu ^{[1]}}
\DN\muN{\mu ^{N}}
\DN\muNc{\check{\mu }^{N}}
\DN\muNx{\mu ^{N}_{x}}
\DN\muNzero{\mu ^{N}_{0}}
\DN\muNone{\mu ^{N,{[1]}}}
\DN\muNm{\mu ^{{N,[m]}}}

\DN\bbb{\mathsf{b}}
\DN\aaa{\mathsf{a}}
\DN\xx{\mathbf{x}}
\DN\xxx{\mathsf{x}}
\DN\yy{\mathbf{y}}
\DN\yyy{\mathsf{y}}
\DN\zz{\mathbf{z}}
\DN\zzz{\mathsf{z}}
\DN\xyi{|x-y_i|}

%Dirichlet form のドメインとか
\DN\di{\mathcal{D} _{\circ}}
\DN\dione{C^{\infty}_{0} (\sS )\otimes \di }
\DN\dik{\di ^{k}}
\DN\dak{\mathcal{D} ^{\aaa ,k}}

\DN\E{\mathcal{E}}
\DN\Ea{\E ^{\aaa }}
\DN\Eak{\E ^{\aaa ,k}}

%log derivative
\DN\dlog{\mathsf{d}}
\DN\dmu{\dlog ^{\mu }}
\DN\dmuN{\dlog ^{N}}

%labeled map
\DN\ulab{\mathfrak{u} }
\DN\lab{\mathfrak{l} }
\DN\lpath{\lab _{\mathrm{path}}}
\DN\lkpath{\lab _{k,\mathrm{path}}}
\DN\upath{\ulab _{\mathrm{path}}}

\DN\PP{\mathbf{P}}
\DN\PPs{\mathbf{P}_{\mathbf{s}}}
\DN\PPP{\mathsf{P}}
\DN\PPPs{\mathsf{P} _{\mathsf{s}}}
\DN\Pmu{\mathsf{P}_{\mu }}
\DN\PmuX{\mathsf{P}_{\mu }^{\XXX _t}}
\DN\Pxt{\PPP _{(x,\sss ) }}
\DN\PPk{\PPP ^{k}}
\DN\Pmt{\PPP _{\sss }}
\DN\Pmul{\PP _{\mu ^{\lab }}}

%log derivativeの有限近似
\DN\g{g}
\DN\gN{\g ^{N}}
\DN\gNrs{\gN _{rs}}
\DN\gNs{\gN _{s}}

\DN\ggN{\mathsf{g} ^{N}}
\DN\ggNrs{\ggrs ^{N}}
\DN\ggNst{\ggst ^{N}}
\DN\ggNrt{\ggrt ^{N}}
\DN\ggNs{\ggN _{s}}
\DN\ggNr{\ggN _{r}}
\DN\ggrs{\mathsf{g}_{rs}} 
\DN\ggst{\mathsf{g}_{st}} 
\DN\ggrt{\mathsf{g}_{rt}} 

\DN\ggr{\mathsf{g}_{r}}
\DN\ggtilder{\ggtilde_{r}}
\DN\ggs{\mathsf{g}_{s}}
\DN\ggtildes{\ggtilde_{s}}
\DN\ggtilde{\tilde{\mathsf{g}}}
\DN\ggtilders{\tilde{\mathsf{g}}_{rs}}

\DN\wNr{{\mathsf w}_r^N}

\DN\hh{\mathfrak{h}}
\DN\hsym{h _{\mathrm{sym}}}

\DN\Y{Y}

\DN\uu{u}
\DN\uN{\uu ^{N}}
\DN\vN{v^{N}}
\DN\vNsinfty{\int_{s\le |x-y|}\vN (x,y)dy }
\DN\vNs{\int_{|x-y|< s}\vN (x,y)dy }

\DN\vvv{v}
\DN\vsinfty{\int_{s \le |x-y|}\vvv (x,y)dy }
\DN\vs{\int_{|x-y|< s}\vvv (x,y)dy }

\DN\rNone{\rho ^{N,\mathrm{1}}}
\DN\rNtwo{\rho ^{N,2}}
\DN\rNnk{\rho ^{N,n+k}}

%FOT関連
\DN\MN{M ^{\nN }}
\DN\sN{\sqrt{\nN }}

\DN\Capa{\mathrm{Cap}}

%Stochastic process

\DN\diai{\diamond i}
 \DN\Xidt{\XXX _{t}^{\diai }}
 \DN\XNidt{\XXX _{t}^{N,\diai }}

\DN\XX{\mathbf{X}}
\DN\XXX{\mathsf{X}}
\DN\XXtt{\mathbf{X}_{\mathbf{t}}}
\DN\Xm{\mathbf{X}^{m}}
\DN\Xms{\mathbf{X}^{m*}}
\DN\XNm{\mathbf{X}^{N,m}}

\DN\SSSz{\SSS _0}
\DN\SSz{\mathbf{S}_0}

\DN\MrTi{\mathsf{M}_{r,T}^i}

\DN\mut{\mu _{\theta }}
\DN\mutx{\mu _{\theta ,x}}
\DN\muNtx{\mu _{\theta ,x}^N}

\DN\muNt{\mu _\theta ^N}
\DN\muNtc{\check{\mu }_\theta ^N}

%SDEgap
\DN\KtwoN{K _{2}^N}
\DN\StN{S _{\theta } ^{N}}

\DN\osc[1]{\psi _{#1}}

\DN\rN{\sqrt{N}}
\DN\orN{\frac{1}{\sqrt{N}}}

\DN\1{\lim_{r \to \infty } \limsup_{N \to \infty} }
\DN\2{\limi{s}\limsupi{N}\sup_{x\in\Sq }}
\DN\3{\frac{\LNyz ^3}{\ellxy \ellxz } }
\DN\4{T_{x,r}^{\nN }}
\DN\5{\frac{1}{\ellxy }}
\DN\8{\frac{1}{|\ellxy |}}
\DN\7{\frac{1}{\ellxyy }}
\DN\6{\frac{1}{\ellxz }}
\DN\0{ \LNyz \LNxyN \LNzxN }

\DN\Trx{T_{r,\infty }^{\nN }(x)}
\DN\Trxx{T_{r,\infty }^{\nN }(x)^2 }

\DN\Srx{S_{r,\infty }(x)}
\DN\Srxx{\Srx ^2 }

\DN\labN{\lab _{\nN }} \DN\labNm{\lab _{\nN ,m}} \DN\labm{\lab _{m }}

\DN\aaaa{\sqrt{2}}
\DN\bbbb{2}

%
                 % Do not remove
%
%
%

%

\begin{abstract}
The distributions of $ N $-particle systems of Gaussian unitary ensembles converge to Sine$ _2$ point processes under bulk-scaling limits. These scalings are parameterized by a macro-position $ \theta $ in the support of the semicircle distribution. The limits are always Sine$ _{2}$ point processes and independent of the macro-position $ \theta $ up to the dilations of determinantal kernels. We prove a dynamical counter part of this fact. 
We prove that the solution of the $ N $-particle systems given by stochastic differential equations (SDEs) converges to the solution of the infinite-dimensional Dyson model. We prove the limit infinite-dimensional SDE (ISDE), referred to as  Dyson's model, is independent of the macro-position $ \theta $, whereas the $ N $-particle SDEs depend on $ \theta $ and are different from the ISDE in the limit whenever $  \theta \not= 0 $. 
\footnote{\textbf{keywords} \texttt{the Gaussian Unitary Ensemble; Dyson's model; bulk scaling limit}}
\end{abstract}

\section{Introduction}\label{s:1}
Gaussian unitary ensembles (GUE) are Gaussian ensembles defined on the space of random matrices $ \MN $ $ (\nN \in \N )$ with independent random variables, the matrices of which are Hermitian. By definition, $ \MN =[\MN _{i,j}]_{i,j=1}^{\nN }$ is then an $ \nN \ts \nN $ matrix having the form 
\begin{align}\notag %\label{:10a}
&
\MN _{i,j} = \begin{cases}
\xi _i&\text{ if } i=j \\ 
\tau  _{i,j}  /\sqrt{2}+ \sqrt{-1} \zeta _{i,j}  /\sqrt{2}
& \text{ if } i<j, 
\end{cases}
\end{align}
where $\{ \xi _{i} ,\tau  _{i,j} , \zeta _{i,j}\}_{i<j}^{\infty }$ 
are i.i.d. Gaussian random variables with mean zero and unit variance. 
Then the eigenvalues $ \lambda_1 ,\ldots, \lambda_{\nN }$ of $ \MN $ 
are real and have distribution $ \muNc $ such that 
\begin{align}\label{:10b}
\muNc (d\xx _N)=\frac{1}{Z^N} \prod_{i<j}^{N}
 | x_i - x_j |^2 \prod_{k=1}^{N} e^{-| x_k |^2 } \,d\xx _N
,\end{align}
where $\xx _N=(x_1,\ldots,x_N)\in \R ^{\nN }$ and 
$Z^N$ is a normalizing constant \cite{AGZ}. 
Wigner's celebrated semicircle law asserts that their empirical distributions converge in distribution to a semicircle distribution:
\begin{align} \notag % \label{:10c}
&
\limi{\nN } \frac{1}{\nN }
\{ \delta_{\lambda_1/\sN } +\cdots + \delta_{\lambda_{\nN /\sN }}  \} 
= \frac{1}{\pi } 1_{(-\aaaa ,\aaaa )}(x) \sqrt{\bbbb -x^2} dx 
.\end{align}
One may regard this convergence as a law of large numbers because the limit distribution is a {\em non-random} probability measure. 

We consider the scaling of the next order in such a way that the distribution is 
supported on the set of configurations. That is, let $ \theta $ be the position of the macro scale given by  
\begin{align}\label{:10d}&
- \aaaa < \theta < \aaaa 
\end{align}
and take the scaling $ x \mapsto y $ such that 
\begin{align}\label{:10e}&
 x = \frac{y}{\sqrt{\nN }} + \theta \sqrt{\nN }
.\end{align}
Let $ \muNt $ be the point process for which the labeled density 
$ \mathbf{m}_{\theta }^{\nN } d\mathbf{x}_{\nN }$ is given by 
\begin{align}\label{:10f}&
\mathbf{m}_{\theta }^{\nN } (\mathbf{x}_{\nN })= 
\frac{1}{Z^N} \prod_{i<j}^{N} | x_i - x_j |^2 \prod_{k=1}^{N} 
e^{-| {x_k} + \theta \nN  |^2/\nN  } 
.\end{align}
The position $ \theta $ in \eqref{:10d} is called the bulk and the scaling in \eqref{:10e} 
the bulk scaling (of the point processes). 
It is well known that the rescaled point processes $  \muNt $ satisfy 
\begin{align}\label{:10h}&
\limi{\nN }  \muNt = \mut \quad \text{ in distribution}
,\end{align}
where $ \mut $ is the determinantal point process with sine kernel $ \mathsf{K}_{\theta }$: 
\begin{align}& \notag %\label{:10i}
\mathsf{K}_{\theta } (x,y) = 
\frac{\sin \{\sqrt{2-\theta ^2} (x-y) \}}{\pi(x-y)}
.\end{align}
By definition $\mu _\theta$ is the point process on $ \R $ 
for which the $ m $-point correlation function $ \rho _{\theta }^m $ 
with respect to the Lebesgue measure is given by 
\begin{align}\notag %\label{:10j}
&
\rho _{\theta }^m (x_1,\ldots,x_m) = 
\det [\mathsf{K}_{\theta } (x_i,x_j) ]_{i,j=1}^m
.\end{align}

We hence see that the limit is universal in the sense that it is 
the Sine$ _{2}$ point process and 
independent of the macro-position $ \theta $ up to the dilation 
of determinantal kernels $ \mathsf{K}_{\theta } $. 
This may be regarded as a first step of the universality of the Sine$ _{2}$ point process, 
which has been extensively studied (see \cite{bey.14-duke} and references therein). 

Once a static universality is established, 
then it is natural to enquire of its dynamical counter part. 
Indeed, we shall prove the dynamical version of \eqref{:10h} and 
present a phenomenon called stochastic-differential-equation (SDE) gaps for $ \theta \not = 0 $. 

Two natural $ \nN $-particle dynamics are known for GUE. 
One is Dyson's Brownian motion corresponding to time-inhomogeneous 
$ \nN $-particle dynamics given by the time evolution of eigenvalues of 
time-dependent Hermitian random matrices $ \mathcal{M}^{\nN }(t)$ 
for which the coefficients are Brownian motions $ B_t^{i,j}$ \cite{mehta}. 

The other is a diffusion process 
$ \mathbf{X}^{\theta ,\nN }=(X^{\theta ,\nN , i})_{i=1}^{\nN }$ 
given by the SDE such that for $ 1 \leq i \leq N $
\begin{align}\label{:10k}
dX_t^{\theta ,N,i} = dB_t^{i} + \sum_{ j \neq i }^{\nN }
 \frac{1}{X_t^{\theta ,N,i} - X_t^{\theta ,N,j}}dt - 
\frac{1}{N}X_t^{\theta ,N,i}\,dt -\theta\, dt 
,\end{align}
which has a unique strong solution for 
$ \mathbf{X}_0^{\theta ,\nN } \in \R ^{\nN } \backslash \mathcal{N}$ 
and $ \mathbf{X}^{\theta ,\nN }$ never hits $ \mathcal{N} $, where 
$ \mathcal{N} = % \cup_{i\not=j}^{\nN }
 \{ \mathbf{x}=(x_k)_{k=1}^N;\, x_i=x_j \text{ for some } i\not=j \} $ \cite{inu}. 

The derivation of \eqref{:10k} is as follows: Let 
$ \muNtc (d\mathbf{x}_{\nN }) 
= \mathbf{m}_{\theta }^{\nN } (\mathbf{x}_{\nN })d\mathbf{x}_{\nN }$ 
 be the labeled symmetric distribution of $ \muNt $. 
Consider a Dirichlet form on $ L^2(\R ^{\nN }, \muNtc )$ such that 
\begin{align}\notag &%\label{:10l}&
\mathcal{E}^{\muNtc } (f,g) = \int_{ \R ^{\nN }} 
\frac{1}{2} \sum_{i=1}^{\nN } \PD{f}{x_i} \PD{g}{x_i} 
\muNtc (d\mathbf{x}_{\nN } )
.\end{align}
Then using \eqref{:10f} and integration by parts, we specify 
the generator $ -A^{\nN }$ of $ \mathcal{E}^{\muNtc }$ on 
$ L^2(\R ^{\nN }, \muNtc )$ such that 
\begin{align}&\notag %\label{:10m}&
A^{\nN } = \half \Delta 
+ \sum_{i=1}^{\nN } \{\sum_{j\not=i }^{\nN }
\frac{1}{x_i-x_j} \} \PD{}{x_i} 
-\sum_{i=1}^{\nN } \{\frac{x_i}{\nN } + \theta \} \PD{}{x_i} 
.\end{align}
From this we deduce that the associated diffusion $ \mathbf{X}^{\theta ,\nN }$ is given by \eqref{:10k}.

Taking the limit $ \nN \to \infty $ in \eqref{:10k}, we {\em intuitively} obtain 
the infinite-dimensional SDE (ISDE) of 
$ \mathbf{X}^{\theta } = (X^{\theta , i})_{i\in\N }$ such that 
\begin{align}\label{:10n}
d X_t^{\theta ,i} &= dB_t^i +\sum_{j\neq i}^{\infty} 
\frac{1}{X_t^{\theta ,i} - X_t^{\theta ,j}}\,dt -\theta \,dt
,\end{align}
which was introduced in \cite{spohn.dyson} with $ \theta = 0 $. 
For each $ \theta $, we have a unique, strong solution 
$ \mathbf{X}^{\theta } $ of \eqref{:10n} such that 
$ \mathbf{X}_0^{\theta } = \mathbf{s} $ 
for $ \mut \circ \lab ^{-1}$-a.s.\! $ \mathbf{s}$, where $ \lab $ is a labeling map. 
Although only the $ \theta = 0 $ ISDE of $ \mathbf{X}^0 =: \mathbf{X} = (X^i)_{i\in\N}$  is studied in \cite{o-t.tail,tsai.14}, 
the general $ \theta \not= 0$ ISDE 
is nevertheless follows easily using the transformation 
\begin{align*}& X_t^{\theta ,i} = X_t^i - \theta t 
.\end{align*}
Let $ \mathsf{X}_t^{\theta } = \sum_i \delta_{X_t^{\theta , i}}$ be the associated delabeled process. Then $ \mathsf{X}^{\theta }= \{ \mathsf{X}_t^{\theta } \} $ takes $ \mut $ as an invariant probability measure, and is {\em not} $ \mut $-symmetric for $ \theta \not= 0$. 

The precise meaning of the drift term in \eqref{:10n} is the substitution of 
$ \mathbf{X}_t^{\theta } =(X_t^{\theta , i})_{i\in\N } $ 
for the function $ b (x,\mathsf{y})$ given by the conditional sum 
\begin{align}\label{:10x}&  b (x,\mathsf{y}) = 
\limi{r} \{ \sum_{|x-y_i|<r} \frac{1}{x-y_i} \} - \theta 
\quad \text{ in 
$ L_{\mathrm{loc}}^1 (\mut ^{[1]})$}
,\end{align}
where $ \mathsf{y}=\sum_i \delta_{y_i}$ and 
$ \mut ^{[1]}$ is the one-Campbell measure of $ \mut $ (see \eqref{:22h}). 
We do this in such a way that 
$ b (X_t^{\theta , i}, \sum_{j\not=i}\delta_{X_t^{\theta , j}})$. 
Because $ \mut $ is translation invariant, 
it can be easily checked that \eqref{:10x} is equivalent to \eqref{:10y}:  
\begin{align}\label{:10y}&  b (x,\mathsf{y}) = 
\limi{r} \{ \sum_{|y_i|<r} \frac{1}{x-y_i} \} - \theta 
\quad \text{ in $ L_{\mathrm{loc}}^1 (\mut ^{[1]})$}
.\end{align}

Let $ \labN $ and $ \lab $ be labeling maps. 
We denote by $ \labNm  $ and $ \labm $
 the first $ m $-components of $ \labN $ and $ \lab $, respectively. 
We assume that, for each $m \in {\mathbb N}$,
\begin{align}\label{:11a}
\lim_{N \to \infty } \mu_\theta ^{N} \circ 
\labNm ^{-1} = \mu_{\theta } \circ \labm ^{-1} \text{ weakly }
.\end{align}
Let $ \mathbf{X}^{\theta ,\nN }=(X^{\theta ,\nN , i})_{i=1}^{\nN }$ 
 and $ \mathbf{X}=(X^i)_{i\in\N }$ 
be solutions of SDEs \eqref{:10k} and \eqref{:10z}, respectively, such that 
\begin{align}\tag{\ref{:10k}}
dX_t^{\theta ,N,i} = \, & dB_t^{i} + \sum_{ j \neq i }^{\nN }
 \frac{1}{X_t^{\theta ,N,i} - X_t^{\theta ,N,j}}dt - 
\frac{1}{N}X_t^{\theta ,N,i}\,dt -\theta\, dt 
%\end{align}\begin{align}
,\\\label{:10z}
dX_t^i = \, & dB_t^i + \limi{r} \sum_{j\neq i,\, |X_t^i - X_t^j|< r }^{\infty} 
\frac{1}{X_t^i - X_t^j}\,dt  %\quad (i\in\N )
.\end{align}
We now state the first main result of the present paper. 
\begin{theorem}\label{l:11} 
Assume \eqref{:10d} and \eqref{:11a}.  
Assume that $ \mathbf{X}_0^{\theta ,\nN } =  \muNt \circ \labN ^{-1} $ in distribution 
and 
$ \mathbf{X}_0 = \mut \circ \lab ^{-1} $  in distribution. 
Then, for each $ m \in \N $, 
\begin{align} \label{:11b}
\lim_{{\it N} \to \infty} (X^{\theta , N ,1},X^{\theta , N ,2},\dots,X^{\theta , N ,m}) =
 (X^{1},X^{2},\ldots,X^{m}) 
\end{align}
weakly in $ C([0,\infty) ,\mathbb{R}^m)$. 
In particular, the limit $ \mathbf{X}=(X^i)_{i\in\N }$ 
does not satisfy \eqref{:10n} for any $ \theta $ other than $ \theta = 0 $. 
\end{theorem}

We next consider non-reversible initial distributions. 
Let  $ \mathbf{X}^N = (X^{N,i})_{i=1}^N$ and 
$ \mathbf{\Y }^{\theta } =(\Y ^{\theta , i})_{i \in \N }$ 
be solutions of \eqref{:11c} and \eqref{:11d}, respectively, such that 
\begin{align}\label{:11c} \quad 
dX_t^{N,i} = \, &dB_t^{i} + \sum_{ j \neq i }^{\nN }
 \frac{1}{X_t^{N,i} - X_t^{N,j}}dt - \frac{1}{N}X_t^{N,i}\,dt 
,\\ \label{:11d}
d \Y _t^{\theta ,i} =  \, &dB_t^i + \limi{r}
\sum_{j\neq i,\, |\Y _t^{\theta ,i} - \Y _t^{\theta ,j}|< r }^{\infty} 
\frac{1}{\Y _t^{\theta ,i} - \Y _t^{\theta ,j}}\,dt + \theta \,dt 
.\end{align}
Note that $  \mathbf{X}^N =  \mathbf{X}^{0 ,\nN }$ and that 
$ \mathbf{X}^N $ is not reversible with respect to 
$ \muNt \circ \labN ^{-1}$ for any $  \theta \not= 0 $. 
We remark that the delabeld process 
$ \mathsf{\Y }^{\theta }= \{\sum_{i\in\N} \delta_{Y_t^{\theta ,i}}\}$ of 
$ \mathbf{\Y }^{\theta }$ 
has invariant probability measure $ \mut $ 
and is {\em not} symmetric with respect to $ \mut $ for $ \theta \not=0 $. 
We state the second main theorem. 
\begin{theorem}\label{l:12} 
Assume \eqref{:10d} and \eqref{:11a}. 
Assume that $  \mathbf{X}_0^N =  \muNt \circ \labN ^{-1} $ in distribution and 
$ \mathbf{\Y }_0^{\theta } = \mut \circ \lab ^{-1} $ in distribution. 
Then for each $ m \in \N $ 
\begin{align} \label{:11f}
\lim_{{\it N} \to \infty} (X^{{\it N},1},X^{{\it N},2},\dots,X^{{\it N},m}) = 
(\Y ^{\theta ,1}, \Y ^{\theta ,2}, \ldots, \Y ^{\theta ,m}) 
\end{align}
weakly in $ C([0,\infty) ,\mathbb{R}^m)$. 
\end{theorem}

%\begin{remark}\label{r:11} 

\begin{itemize}
\item 
We refer to the second claim in \tref{l:11}, and \eqref{:11f} as the SDE gaps. 
The convergence in \eqref{:11f} of \tref{l:12} resembles the 
\lq\lq Propagation of Chaos'' in the sense that the limit equation \eqref{:11d} 
depends on the initial distribution, although it is a linear equation. 
Because the logarithmic potential is by its nature long-ranged, the effect 
of initial distributions $ \muNt $ still remains in the limit ISDE, and the rigidity of the
Sine$ _{\mathrm{2}}$ point process makes the residual effect a non-random drift term $ \theta dt $. 
%\end{document}
\item 
Let $ \SSS_{\theta}$ be a Borel set such that $ \mut (\SSS_{\theta})= 1 $, 
where $ -\aaaa < \theta < \aaaa $. 
In \cite{kawa.dp}, the first author proves that one can choose $ \SSS_{\theta}$ such that 
$ \SSS _{\theta } \cap \SSS _{\theta '} = \emptyset $ 
if $ \theta \not = \theta '$ and that for each $ \mathsf{s} \in \SSS _{\theta }$ 
\eqref{:10z} has a strong solution $ \mathbf{X}$ such that 
$ \mathbf{X}= \lab (\mathsf{s})$ and that 
%
%\begin{align}\notag %\label{:12a}
%& \quad \quad 
$$\mathsf{X}_t := \sum_{i=1}^{\infty} \delta_{X_t^i} \in \mathsf{S}_{\theta } 
\quad \text{ for all }t \in [0,\infty )
.$$%\end{align}
This implies that the state space of solutions of \eqref{:10z} can be decomposed into 
uncountable disjoint components. We conjecture that the component $ \SSS _{\theta }$ 
is ergodic for each $ \theta \in (-\aaaa , \aaaa )$. 
\item 
For $ \theta = 0 $, 
the convergence \eqref{:11b} is also proved in \cite{o-t.sm}. The proof in \cite{o-t.sm} is algebraic and valid only for dimension $ d=1 $ and inverse temperature $ \beta = 2$ with the logarithmic potential. It relies on an explicit calculation of the space-time correlation functions, the strong Markov property of the stochastic dynamics given by the algebraic construction, the identity of the associated Dirichlet forms constructed by two completely different methods, and the uniqueness of solutions of ISDE \eqref{:10n}. 
\\
Although one may prove \eqref{:11a} for $ \theta \not = 0$ using the algebraic method in \cite{o-t.sm}, this requires a lot of work as mentioned above. We remark that, as a corollary and an application,  \tref{l:11} proves the weak convergence of finite-dimensional distributions explicitly given by the  space-time correlation functions. We refer to \cite{KT10a,o-t.sm} for the representation of these correlation functions.  
\item 
Tsai proves the pathwise uniqueness and the existence of strong solutions of 
\begin{align}\label{:12b}
dX_t^i &= dB_t^i + 
\frac{\beta }{2}\limi{r} \sum_{j\neq i,\, |X_t^i - X_t^j|< r }^{\infty} 
\frac{1}{X_t^i - X_t^j}\,dt \quad (i\in\N )
\end{align}
for general $ \beta \in [1,\infty )$ in \cite{tsai.14}. 
The proof uses the classical stochastic analysis and crucially depends on a specific monotonicity of SDEs \eqref{:12b}. 
For $ \beta = 1,4$, we have a good control of the correlation functions as for $ \beta = 2 $. Hence our method can be applied to $ \beta = 1,4 $ and the same result as for $ \beta = 2 $ in \tref{l:11} holds. We shall return to this point in future. 
%\end{remark}
\end{itemize}

The key point of the proof of \tref{l:11} is to prove the convergence of the drift coefficient $ b ^N(x,\mathsf{y})$ of the $ N $-particle system to the drift coefficient  
$ b (x,\mathsf{y})$ of the limit ISDE even if $ \theta \not= 0 $. That is, as $ N \to \infty $, 
\begin{align}& \notag %\label{:12c}&
b ^N(x,\mathsf{y})= 
\{ \sum_{i=1}^N \frac{1}{x-y_i} \} - \theta \quad 
\Longrightarrow \quad  b (x,\mathsf{y})= 
\limi{r} \{ \sum_{|y_i|<r} \frac{1}{x-y_i} \} 
.\end{align}
Note that support of the coefficients $ b ^N (x,\mathsf{y})$ and $ b (x,\mathsf{y})$ are mutually disjoint, and that the sum in $ b ^N $ is not neutral for any $ \theta \not= 0 $. We shall prove uniform bounds of the tail of the coefficients using fine estimates of the correlation functions, and cancel out the deviation of the sum in $ b ^N $ with $ \theta $. 
Because of rigidity of the Sine$_{\mathrm{2}}$ point process, we justify this cancellation not only for static but also dynamical instances. 

The organization of the paper is as follows: In \sref{s:2}, we prepare general theories for interacting Brownian motion in infinite dimensions. In \sref{s:3}, we quote estimates for the oscillator wave functions and determinantal kernels. In \sref{s:4}, we prove key estimates \eqref{:26a}--\eqref{:26d}. In \sref{s:5}, we complete the proof of \tref{l:11}.  
In \sref{s:6}, we prove \tref{l:12}.

\section{Preliminaries from general theory}\label{s:2}
In this section we present the general theory described in \cite{o.isde,o.rm,o-t.tail,k-o.fpa} 
in a reduced form sufficient for the current purpose. In particular, we take the space where particles move in $ \R $ rather than $ \Rd $ as in the cited articles. 

\subsection{$ \mu $-reversible diffusions} \label{s:21}
Let $\Sr =\{s \in \R \, ;  |s| < r  \} $. 
The configuration space $\SSS $ over $ \R $ is a Polish space equipped with the vague topology such that 
\begin{align*}&
 \SSS  = \{ \mathsf{s} = \sum _i \delta _{s_i}\, ;\, \mathsf{s} ( \Sr ) < \infty 
\text{ for all } r \in \N \} 
.\end{align*}
Each element $ \mathsf{s} \in \SSS $ is called a configuration regarded as countable delabeled particles. A probability measure $\mu $ on 
$(\SSS ,\mathcal{B}(\SSS ) )$ is called a point process (a random point field). 

A locally integrable symmetric function 
$\map{\rho ^n}{\sS  ^n}{[0 ,\infty)}$ is called the $ n $-point correlation function of $\mu $ with respect to the Lebesgue measure if $\rho ^n$ satisfies

\begin{align}\notag % \label{:20a}
\int_{A_1^{k_1}\times \cdots \times A_m^{k_m}} 
\rho ^n (s_1,\ldots,s_n)  \,d\mathbf{s}_n 
 = \int _{\SSS } \prod _{i = 1}^{m} \frac{\mathsf{s} (A_i) ! }
{(\mathsf{s} (A_i) - k_i )!} \mu (d\mathsf{s})
\end{align}
for any sequence of disjoint bounded measurable subsets 
$ A_1,\ldots,A_m \subset \sS  $ and a sequence of natural numbers 
$ k_1,\ldots,k_m $ satisfying $ k_1+\cdots + k_m = n $. Here we assume that 
${\mathsf{s} (A_i) ! }/{(\mathsf{s} (A_i) - k_i )!} =0$ for $\mathsf{s} (A_i) - k_i < 0$. 

Let $ \map{\Phi }{\sS  }{\R }$ and $ \map{\Psi }{\sS  ^2}{\R \cup \{ \infty \} }$ be 
measurable functions called free and interaction potentials, respectively. 
Let $\mathcal{H}_r$ be the Hamiltonian on $\Sr $ given by
\begin{align}\notag %\label{:20c}
\mathcal{H}_r(\xxx ) =\sum_{x_i\in\Sr }\Phi (x_i) +\sum_{j\neq k,x_j,x_k\in\Sr }\Psi (x_j,x_k)\quad \text{ for }\xxx =\sum_i \delta_{x_i }
.\end{align}
For each $ m ,r \in \N $ and $ \mu $-a.s.\! $ \xi \in \SSS $, let  
$ \mu _{r,\xi }^{m} $ denote the regular conditional probability such that 
\begin{align}\notag 
\mu _{r,\xi }^{m} =\mu (\pi _{\Sr }\in\cdot \,|\, \pi _{\Sr ^c}(\xxx )=\pi _{\Sr ^c}(\xi ),\, \xxx (\Sr )=m)
.\end{align}
Here for a subset $ A $, we set $ \map{\pi_{A}}{\SSS }{\SSS } $ by 
 $ \pi_{A} (\mathsf{s}) = \mathsf{s} (\cdot \cap A )$.

Let $\Lambda _r$ denote the Poisson point process with intensity 
being a Lebesgue measure on $\Sr $. 
We set  $\Lambda _r^m (\cdot)=\Lambda _r(\cdot \cap \Srm)$, where 
$\Srm =\{\sss \in\SSS \,;\, \sss (\Sr )=m \}$. 
\begin{definition}[\cite{o.rm}, \cite{o.rm2}] \label{d:21}
A point process $\mu $ is said to be a $(\Phi ,\Psi )$-quasi Gibbs measure if 
its regular conditional probabilities $ \mu _{r,\xi }^{m} $ satisfy,  
 for any $r,m\in\N$ and $\mu $-a.s. $\xi $, 
\begin{align}\notag %\label{:20e}
\cref{;22a}^{-1} e^{-\mathcal{H}_r(\xxx )}\Lambda _r^m (d\xxx ) \le\mu _{r,\xi }^{m}(d\xxx ) \le \cref{;22a} e^{-\mathcal{H}_r(\xxx )}\Lambda _r^m (d\xxx )
.\end{align}
Here $\Ct{;22a}$ is a positive constant depending on $r,m,\xi $. 
\end{definition}

The significance of the quasi-Gibbs property is to guarantee the existence of $ \mu $-reversible diffusion process $ \{ P_{\mathsf{s}} \} $ on $ \SSS $ 
given by the natural Dirichlet form associated with $ \mu $, 
in analogy with distorted Brownian motion in finite-dimensions. 

To introduce the Dirichlet form, we provide some notations. 
 We say a function $ f $ on $ \SSS $ is local if 
$ f $ is $ \sigma [\pi_{K}]$-measurable for some compact set $ K $ in $ \sS  $. 
 For a local function $ f $ on $ \SSS $, we say $ f $ is smooth if $ \check{f}$ is smooth, where 
 $ \check{f}(x_1,\ldots )$ is the symmetric function such that 
 $ \check{f}(x_1,\ldots ) = f (\mathsf{x})$ for $ \mathsf{x} = \sum _i \delta _{x_i}$. 
 Let $ \di $ be the set of all bounded, locally smooth functions on $ \SSS $.

Let $ \mathbb{D}$ be the standard square field on $ \SSS $ such that 
for $ f , g \in \di $ and $ \mathsf{s}=\sum_i\delta_{s_i}$
\begin{align*}&
\mathbb{D}[f,g] (\mathsf{s})= 
\half \{\sum_i (\nabla _i{\check{f}} )(\nabla _i{\check{g}}) \} \, 
 (\mathsf{s})
.\end{align*}
We write $ \mathbf{s}=(s_i)_i$. Because the function 
$ \sum_i (\nabla _i {\check{f}}) (\mathbf{s}) (\nabla _i{\check{g}})  (\mathbf{s})$ 
is symmetric in $ \mathbf{s}=(s_i)_i$, we regard it as a function of $ \mathsf{s}$. 
We set $ \Lm = L^2(\SSS ,\mu )$ and let 
\begin{align*}&
\mathcal{E}^{\mu }(f,g) = \int_{\SSS } \mathbb{D}[f,g] (\mathsf{s}) \mu (d\mathsf{s}), \quad 
%\\ & 
 \di ^{\mu } =\{ f \in \di \cap \Lm \, ;\, \mathcal{E}^{\mu }(f,f) < \infty \} 
.\end{align*}
We quote: 
\begin{lemma}[{\cite{o.rm}}]	\label{l:21}
Assume that $ \mu $ is a $(\Phi ,\Psi )$-quasi Gibbs measure 
with upper semicontinuous $(\Phi ,\Psi )$. Assume that the correlation functions $ \{ \rho^n \} $ are locally bounded for all $ n \in \N $. 
Then $ (\mathcal{E}^{\mu }, \di ^{\mu } )$ is closable on $ \Lm $. 
Furthermore, there exists a 
$ \mu $-reversible diffusion process $ \{ P_{\mathsf{s}} \} $ associate with the Dirichlet form $ (\mathcal{E}^{\mu }, \mathcal{D}^{\mu } )$ on $ \Lm $. 
Here $ (\mathcal{E}^{\mu }, \mathcal{D}^{\mu } )$ is the closure of 
$ (\mathcal{E}^{\mu }, \di ^{\mu })$ on $ \Lm $. 
\end{lemma}

\subsection{Infinite-dimensional SDEs}
\label{s:22}

Suppose that diffusion $ \{ P_{\mathsf{s}} \} $ in \lref{l:21} 
is collision-free and that each tagged particle does not explode. 
Then we can construct labeled dynamics $ \mathbf{X}=(X^i)_{i\in\Z }$ 
by introducing the initial labeling $ \lab =(\lab _i)_{i\in\Z }$ such that 
\begin{align}\notag %\label{:22p}
&
\mathbf{X}_0 = \lab (\mathsf{X}_0)
.\end{align}
Indeed, once the label $ \lab $ is given at time zero, then each particle retains the tag for 
all time because of the collision-free and explosion-free property. 

To specify the ISDEs satisfied by $ \mathbf{X}$ above, 
we introduce the notion of the logarithmic derivative of 
$ \mu $, which was introduced in \cite{o.isde}. 

A point process $ \mu _{x}$ is called the reduced Palm measure of $ \mu $ 
conditioned at $ x \in \sS  $ if $ \mu _{x}$ is the regular conditional probability defined as 
\begin{align} & \notag %\label{:22j}
  \mu _{x} = \mu (\cdot - \delta_x | \mathsf{s} (\{ x \} ) \ge 1 )
.\end{align} 
A Radon measure $ \muone $ on $ \sS  \times \SSS $ 
is called the 1-Campbell measure of $ \mu $ if 
\begin{align}\label{:22h}&
 \muone  (dx d\mathsf{s}) = \rho ^1 (x) \mu _{x} (d\mathsf{s}) dx 
.\end{align}
We write $ f \in L_{\mathrm{loc}}^p(\muone )$ if 
$ f \in L^p(\Sr \ts \SSS , \muone  )$ for all $ r \in\N $. 
\begin{definition}\label{d:22}
A $ \R $-valued function  $ \dmu \in L_{\mathrm{loc}}^1(\muone  )$  is called the 
 {\em logarithmic derivative} of $\mu $ if, for all 
$\varphi \in C_{0}^{\infty}(\sS  )\otimes \di $, 
\begin{align}&\notag % \label{:21i}&
 \int _{\sS  \times \SSS } 
 \dmu   (x,\mathsf{y})\varphi (x,\mathsf{y})   \muone  (dx d\mathsf{y})  =  
  - \int _{\sS  \times \SSS }  \nabla  _x \varphi (x,\mathsf{y})  \muone  (dx d\mathsf{y}) 
.\end{align}
\end{definition}
Under these assumptions, we obtain the following: 
\begin{lemma}[{\cite{o.isde}}]	\label{l:22}
Assume that $ \mathbf{X}=(X^i)_{i\in\N }$ is the collision-free and explosion-free.  
Then $ \mathbf{X}$ is a solution of the following ISDE: 
\begin{align}\label{:22a}&
\quad dX_t^i = dB_t^i + \frac{1}{2} \dmu (X_t^i , \Xidt ) dt \quad (i \in\N )
\end{align}
with initial condition $ \mathbf{X}_0 = \mathbf{s}$ for $ \mul $-a.s.\! $ \mathbf{s}$, where $ \Xidt = \sum_{j\not=i}^{\infty} \delta_{X_t^{j}}$. 
\end{lemma}

\subsection{Finite-particle approximations}\label{s:23} 
Let $ \mu $ be a point process with correlaton functions $ \{ \rho ^n  \}_{n\in\N } $. 
Let $ \{\muN \}_{N \in \N } $ be a sequence of point processes on $\sS  $ 
such that $ \muN (\{ \mathsf{s}(\sS  ) = \nN \} ) = 1 $.  We assume: \\
\As{A1} Each $ \muN $ has correlation functions $ \{\rho ^{N,n}\}_{n \in \N } $ 
satisfying, for each $r \in\N$, 
\begin{align}\label{:23b} & 
\limi{N} \rho ^{N,n} (\mathbf{x})= 
\rho ^{n} (\mathbf{x}) \quad \text{ uniformly on $\Sr ^{n}$ for each $n\in\N$}  
,\\ \label{:23c}& 
\sup_{N\in\N } \sup_{\mathbf{x} \in \Sr ^{n}} \rho ^{N,n} (\mathbf{x}) \le 
\cref{;40b} ^{n} n ^{\cref{;40c}n} 
,\end{align}
where $ 0 < \Ct{;40b}(r) < \infty $ and 
$ 0 < \Ct{;40c}(r)< 1 $ are constants independent of $ n \in \N $. 

\medskip 

It is known that \eqref{:23b} and \eqref{:23c} imply the weak convergence 
of $ \{ \muN  \} $ to $ \mu $ \cite[Lemma A.1]{o.rm}. As in \sref{s:1}, 
let $ \lab $ and $ \labN $ be  labels of $ \mu $ and $ \muN $, respectively. 
We assume: 

\medskip
\noindent
\As{A2} For each $m\in\N$, 
 \begin{align}\notag %\label{:23d}
& \quad 
 \limi{N}\mu^{\nN } \circ \labNm ^{-1} =\mu \circ \labm ^{-1}
 \quad \text{ weakly in $ \sS  ^m $}
 .\end{align}

\medskip 

We shall later take $ \mu^{\nN } \circ \labN ^{-1} $ 
as an initial distribution of labeled finite particle system. 
Therefore, \As{A2} means the convergence of the initial distribution of the labeled dynamics. 

For a labeled process $ \mathbf{X}^N=(X^{N,i})_{i=1}^{\nN }$, where 
$ \nN \in \N $, we set  
\begin{align} & \notag %\label{:20l}&
 \mathsf{X}_t^{N,\diai }= \sum_{j\not=i}^{ N } \delta_{X_t^{N,j}} 
,\end{align}
where $ \mathsf{X}_t^{N,\diai }$ denotes the zero measure for $ \nN = 1 $. 
Let $ \map{\bbb ^{\nN } ,\bbb }{\sS  \ts \SSS }{\R }$ be measurable functions. 
We introduce the finite-dimensional  SDE 
of $ \mathbf{X}^{\nN } =(X^{\nN ,i})_{i=1}^{\nN } $ 
with these coefficients such that for $ 1\le i\le N $ 
\begin{align}\label{:23f}
dX_t^{N,i} &= dB_t^i  + \bbb ^{\nN }(X_t^{N,i},\XNidt )dt 
.\end{align}
We assume: 

\medskip

\noindent \As{A3} 
 SDE \eqref{:23f} with initial condition 
$ \XX _0^{\nN }= \mathbf{s}$
has a unique solution for $ \muN \circ \labN ^{-1}$-a.s.\! $ \mathbf{s}$ 
for each $ \nN $. 
This solution does not explode. 

\medskip

Let $ \map{\uu ,\ \uN ,\ \ww }{\sS  }{\R }$ and $\map{\g }{\sS  ^{2} }{\R }$ 
 be measurable functions. We set 
\begin{align}
\label{:23g}&
\ggr (x,\mathsf{y})= \sum_{i} \chi _r  (x-y_i) \g (x,y_i)
,\\&\label{:23h}
\rrNs (x,\mathsf{y})= \sum_{i}(1- \chi _r  (x-y_i))\g (x,y_i) 
,\end{align}
where $ \mathsf{y}=\sum_{i}\delta_{y_i}$ and 
$ \chi _r  \in C_0^{\infty} (\sS  )$ is a cut-off function such that 
$ 0 \le \chi _r  \le 1 $, $ \chi _r  (x) = 0 $ for $ |x| \ge r + 1 $, and 
$ \chi _r  (x) = 1$ for $ |x| \le r $. 
We assume the following. 

\medskip

\noindent 
\As{A4} Each $ \muN $ has a logarithmic derivative $ \dmuN $ such that 
\begin{align}\label{:23i}&
\dmuN (x,\mathsf{y})= \uN (x) + \ggr (x,\mathsf{y}) + \rrNs (x,\mathsf{y})
.\end{align}
Furthermore, we assume that 
\begin{itemize}

\item[(1)]
$\uN $ are in $ C^1 (\sS  ) $. Furthermore, 
$ \uN $ and $ \nabla \uN $ converge uniformly to $ u $ and $ \nabla u $, respectively, on each compact set in $ \sS  $.

\item[(2)] 
$ g \in C^1 (\sS  ^2 \cap \{x \not= y \}) $. 
There exists a $ \phat >1 $ such that, for each $ R \in \N $, 
\begin{align}\label{:23j}&
\limi{\p } \limsupi{\nN } 
\int_{x \in \SR , |x-y| \le 2^{-\p } } 
\chi _r  (x-y) | g (x,y) | ^{\phat }\, \rho _x^{\nN ,1}(y) dxdy = 0 
,\end{align}
where $ \rho _x^{\nN ,1} $ is a one-correlation function of the reduced Palm measure 
$ \mu _x^{\nN }$. 

\item[(3)] There exists a continuous function $ \map{w }{\sS  }{\R }$ such that 
for each $ R \in \N $
\begin{align}\label{:23k}& 
\limi{r}\limsupi{N} 
 \int_{\SR \times \SSS } |\rrNs (x,\mathsf{y}) - \www |^{\phat } d\muNone 
= 0 % , \quad \ww \in L^{\phat }_{\mathrm{loc}} (\sS  ,dx)
.\end{align}
\end{itemize}

Let $ p $ be such that $ 1 < p < \phat $. 
Assume \As{A1} and \As{A4}. 
Then \cite[Theorem 45]{o.isde} deduces that 
the logarithmic derivative $\dmu $ of $ \mu $ exists in 
$ \Lploc (\muone )$ and is given by 
\begin{align}& \label{:23l}
\dmu (x,\mathsf{y})= 
\uu (x) + \ff  (x,\mathsf{y}) + \www 
.\end{align}
Here $ \ff  (x,\mathsf{y})= \limi{r} \ggr (x,\mathsf{y}) $ and 
the convergence of $ \lim \ggr $ takes place in $ \Lploc (\muone )$. 
Taking \eqref{:23l} into account, 
we introduce the ISDE of $ \mathbf{X}=(X^i)_{i\in\N }$: 
\begin{align}\label{:23m}
dX_t^i &= dB_t^i + \frac{1}{2} \{ u (X_t^i) + 
\ff  (X_t^i,\mathsf{X}_t^{\diai }) + w (X_t^i) \} dt 
.\end{align}

Under the assumptions of \lref{l:22}, ISDE \eqref{:23m} 
with $ \mathbf{X}_0 = \mathbf{s}$ 
has a solution for $ \mul $-a.s.\!\! $ \mathbf{s}$. Moreover, the associated delabeled diffusion $ \mathsf{X}=\{ \mathsf{X}_t \} $ is $ \mu $-reversible, where 
$ \mathsf{X}_t = \sum_{i\in\N } \delta_{X_t^i}$ 
for $ \mathbf{X}_t = (X_t^i)_{i\in\N }$. 
As for uniqueness, we recall the notion of $ \mu $-absolute continuity solution 
introduced in \cite{o-t.tail}. 

Let $ \mathbf{X} = (X^i)_{i\in\N }$ be a family of solution of \eqref{:23m} satisfying 
$ \mathbf{X}_0 = \mathbf{s} $ for $ \mu \circ \lab ^{-1}$-a.s.\! $ \mathbf{s}$. 
Let $ \mu _t $ be the distribution of the delabeled process 
$ \mathsf{X}_t =\sum_{i\in\N } \delta_{X_t^i}$ at time $ t $ 
with initial distribution $ \mu $. That is, $ \mu _t $ is given by 
\begin{align*}&
\mu _t = \int_{\SSS } P_{\mathsf{s}} (\mathsf{X}_t \in \cdot ) d\mu 
\end{align*}
We say that $ \mathbf{X}$ satisfies the $ \mu $-absolute continuity condition if 
\begin{align}\label{:23o}&
\mu _t \prec \mu \quad \text{ for all } t \ge 0 
,\end{align}
where $ \mu _t \prec \mu $ means that $ \mu _t $ 
is absolutely continuous with respect to $ \mu $. 
If $ \mathsf{X}$ is $ \mu $-reversible, then \eqref{:23o} is satisfied. 

We say ISDE \eqref{:23m} has $ \mu $-uniqueness in law of solutions 
if $ \mathbf{X}$ and $ \mathbf{X}'$ are solutions with the same initial distributions 
 satisfying the $ \mu $-absolute continuity condition, then 
 they are equivalent in law. 
We assume: 

\medskip 

\noindent 
\As{A5} 
 ISDE \eqref{:23m} has $ \mu $-uniqueness in law of solutions.

\medskip

It is proved in \cite{o-t.tail} that ISDE \eqref{:22a} has a $ \mu $-pathwise unique strong solution if $ \mu $  is tail trivial, the logarithmic derivative $ \dmu $ has a sort of off-diagonal smoothness, and the one-correlation function has sub-exponential growth at infinity. 
This results  implies $ \mu $-uniqueness in law. 
We refer to Theorems 2.1 and 9.3 in \cite{o-t.tail} for details. 
The next result is a special case of \cite[{Theorem 2.1}]{k-o.fpa}. 
\begin{lemma}[{\cite[{Theorem 2.1}]{k-o.fpa}}]\label{l:23} 
Make the same assumptions in \lref{l:21} and \lref{l:22}. 
Assume \As{A1}--\As{A4}. 
Assume that $ \mathbf{X}_0^{\nN } = \muN \circ \labN ^{-1}$ in distribution. 
Then $ \{ \mathbf{X}^{\nN } \}_{N \in \N } $ is tight in $ C([0,\infty );\sS  ^{\N })$ and 
each limit point $ \mathbf{X}$ of $ \{ \mathbf{X}^{N} \}_{N \in \N } $ is a solution of 
\eqref{:23m} with initial distribution $ \mul $. 
If, in addition, we assume \As{A5}, then for any $m \in \mathbb{N}$ 
\begin{align}\notag %\label{:23a}
& \limi{\nN } %\mathbf{X}^{\nN , m } = \mathbf{X}^m 
(X^{\nN ,1},\ldots,X^{\nN ,m}) = (X^1,\ldots,X^m)
.\end{align}
 weakly in $ C([0,\infty) ,\sS  ^m)$. 
Here $ \mathbf{X}^N=(X^{N,i})_{i=1}^{\nN }$ and 
$ \mathbf{X} = (X^i)_{i\in\N }$ as before. 
\end{lemma}

\subsection{Reduction of \tref{l:11} to \eqref{:23k}} \label{s:24}

In this subsection, we deduce \tref{l:11} from \lref{l:23} by assuming \eqref{:23k}. 
We take $\muNt $ and $\mut $ as in \sref{s:1}. 
Then the logarithmic derivative $  \dlog ^{\muNt }$ of $\muNt $ is given by
\begin{align} \label{:24a}
 \dlog ^{\muNt }(x,\yyy )= \sum_{i=1}^N \frac{2}{x - y_i} -\frac{2x}{N} -2\theta 
,\end{align}
where $ \yyy =\sum_{i} \delta _{y_i}$. 
From \eqref{:24a}, we take coefficients in \As{A4} as follows: 
\begin{align} \label{:24b} 
\uN (x)&= -\frac{2x}{N} -2\theta , \quad  u (x) = - 2 \theta  ,\quad w (x) = 2\theta 
,\\ \label{:24c} 
\g (x,y) & =\frac{2}{x - y} 
%,\\ \label{:24d} w(x) &=2\theta 
.\end{align}
Other functions are given by \eqref{:23g} and \eqref{:23h}. 
 
\begin{lemma} \label{l:24}
Assume \eqref{:23k} holds with $ \phat = 2$ 
for the coefficients as above. Then \eqref{:11b} holds. 
\end{lemma}
\PROOF 
To prove \lref{l:24}, we check the assumptions in \lref{l:23}, that is, 
the assumptions in \lref{l:21}, \lref{l:22}, and \As{A1}--\As{A5}.  

The assumptions in \lref{l:21} are proved in \cite{o.rm}. 
The assumptions in \lref{l:22} are checked in \cite{o.isde}. 
\As{A1} is well known. \As{A2} is assumed by \eqref{:11a}. 
\As{A3} is obvious as the interaction is smooth outside the origin, and 
the capacity of the colliding set $ \{ x_i=x_j \text{ for some }i\not=j\} $ 
is zero (see \cite{o.col,inu}). Furthermore, the one-correlation functions are bounded, 
which guarantees explosion-free of tagged particles. 
We take functions in \As{A4} as \eqref{:24b} and \eqref{:24c}.  
These satisfy \eqref{:23i}, \eqref{:23j}, and \thetag{1} of \As{A4}.  
\eqref{:23k} is satisfied by assumption. It is known that $ \mut $ is tail trivial \cite{o-o.tt}. 
Then \As{A5} follows from tail triviality of $ \mut $ and \cite[Theorem 3.1]{o-t.tail}. 
All the assumptions in \lref{l:23} are thus satisfied, and hence yields \eqref{:11b}.  
\PROOFEND

\subsection{A sufficient condition for \eqref{:23k}} \label{s:25}
The most crucial step to apply \lref{l:23} is to check \eqref{:23k}. 
Indeed, it only remains to prove \eqref{:23k} for \tref{l:11}. 
We quote then a sufficient condition for \eqref{:23k} 
in terms of correlation functions from \cite{o.isde}. 
\lref{l:26} below is a special case of \cite[Lemma 53]{o.isde}. 

Let $ \muNtx $ be the reduced Palm measure of $ \muNt $ 
conditioned at $  x $. 
We denote the supremum norm in $ x $ over $ S _R $ by $ \|\, {\cdot }\, \|_R $. 
Let $ \mathrm{E}^{\cdot }$ and $ \mathrm{Var}^{\cdot }$ denote the expectation and variance with resoect to $ \cdot $, respectively. %
\begin{lemma} \label{l:25} 
Assume $ | \theta | < \sqrt{2} $. Let $ w_r $ be as in \eqref{:23h} 
with $ g (x,y)$ given by \eqref{:24c}. 
Let $ w (x) = 2\theta $ as in \eqref{:24b}.  
Then \eqref{:23k} follows from \eqref{:25a}--\eqref{:25d}.   
\begin{align} 
& \label{:25a}
\1  \SSup{ \mathrm{E}^{\muNt }[w_r(x,\mathsf{y})] - 2\theta  }  =0
,\\& \label{:25b}
\1 \SSup{ \mathrm{E}^{\muNt } [w_r(x,\mathsf{y})] -
\mathrm{E}^{\muNtx } [w_r(x,\mathsf{y})] }  =0
,\\& \label{:25c}
\1  \SSup{ \mathrm{Var}^{\muNt }[w_r(x,\mathsf{y})] }  = 0
,\\&\label{:25d}
\1 \SSup{ \mathrm{Var}^{\muNt } [w_r(x,\mathsf{y})] -
\mathrm{Var}^{\muNtx } [w_r(x,\mathsf{y})] }  = 0
.\end{align}
\end{lemma}
\PROOF 
\lref{l:25} follows from \cite[Lemma 52]{o.isde}. 
Indeed, \eqref{:25a}, \eqref{:25b}, \eqref{:25c}, and \eqref{:25d} in the present paper correspond to \thetag{5.4}, \thetag{5.2}, \thetag{5.5}, and \thetag{5.3} 
 in \cite{o.isde}, respectively. We note that in \cite{o.isde} we use 
$ 1_{\Sr }(x)$  instead of $ \chi _r (x) $. This slight modification yields no difficulty. 
\PROOFEND 

Multiplying $ w_r (x,\mathsf{y})$ by a half, we give a sufficient condition of \eqref{:25a}--\eqref{:25d} in terms of correlation functions. 
Let $ \rho_{\theta ,x}^{N,m}$ and $ \rho_{\theta }^{N,m}$ 
 be the $ m $-point correlation functions of $ \muNtx $ and $ \muNt $, respectively. Let 
\begin{align}\notag %\label{:26z}
&
\Srx = \{y \in \R \, ; r < |x-y| < \infty  \} 
.\end{align}

\begin{lemma} \label{l:26} 
Assume $ | \theta | < \sqrt{2} $. Then \eqref{:25a}--\eqref{:25d} follow from 
\eqref{:26a}--\eqref{:26d}.   
\begin{align} 
& \label{:26a}
\1  \SSup{\int_{\Srx } \frac{\rho_\theta ^{N,1}(y)}{x-y} dy - \theta  }  =0
,\\& \label{:26b}
\1 \SSup{ \int_{\Srx } 
\frac{\rho_{\theta ,x}^{N,1}(y) - \rho_\theta ^{N,1}(y) }{x-y}  \,dy }  =0
,\\& \label{:26c}
\1 \SSup{ \int_{\Srx }  \frac{\rho_\theta ^{N,1}(y)}{(x-y)^2}  dy
 + 
\int_{\Srxx } 
\frac{\rho_{\theta }^{N,2}(y,z)-\rho_\theta ^{N,1}(y) \rho_\theta ^{N,1} (z) }
{(x-y)(x-z)}  \,dy dz }
  = 0
,\\&\label{:26d}
\1 \SSup{
\int_{\Srx } \frac{\rho_{\theta ,x}^{N,1}(y) -\rho_\theta ^{N,1}(y) }{(x-y)^2} 
 \,dy 
\\& \notag 
+\int_{\Srxx } \frac{\rho_{\theta ,x}^{N,2}(y,z)-\rho_{\theta ,x} ^{N,1}(y) \rho_{\theta ,x} ^{N,1} (z)
- \{\rho_\theta ^{N,2}(y,z) -\rho_\theta ^{N,1}(y) \rho_\theta ^{N,1} (z) \} }
{(x-y)(x-z)} 
\,dy dz  }  =0
.\end{align}
\end{lemma}
\PROOF 
\lref{l:26} follows immediately from a standard calculation of correlation functions and the definitions of $ w_r$ and $  \chi _r $. 
\PROOFEND

\section{Subsidiary estimates} \label{s:3}

Keeping \lref{l:26} in mind, our task is to prove \eqref{:26a}--\eqref{:26d}. 
To control the correlation functions in \lref{l:26} we prepare in this section estimates of the oscillator wave functions and determinantal kernels.  
We shall use these estimates in \sref{s:4}. 

\subsection{Oscillator wave functions} \label{s:31}

Let $H_n(x)=(-1)^n e^{x^2} (\frac{d}{dx})^n e^{-x^2} $ be Hermite polynomials.
Let $\osc n(x)$ denote the oscillator wave functions defined by 
\begin{align}\notag %\label{:30a}
&
\osc n(x) =\frac{1}{\sqrt{\sqrt{\pi}2^{n}n!}}e^{-\frac{x^2}{2}}H_n(x)
.\end{align}
Note that $ \{ \osc {n} \}_{n=0}^{\infty} $ is an orthonormal system; 
$\int _{\R} \osc {n}(x)\osc {m}(x)\,dx=\delta _{nm}$. 

The following estimates for these oscillator wave functions 
are essentially due to Plancherel-Rotach \cite{PR29}. 
We quote here a version from Katori-Tanemura \cite{KT11}. 
\begin{lemma}[\cite{KT11}]\label{l:31} 
Let $ C_{nm}^1 $, $ C_{nm}^2 $, and $ D_{nm}^1$ be the constants 
introduced in \cite{KT11} (see \As{A.1} in \cite[572 p]{KT11}). 
Let $ l = -1,0,1$ and $N , L \in \mathbb{N} $. Then we have the following. \\
\thetag{1} Let $0 < \tau \le \frac{\pi}{2} $. 
Assume that $N\sin ^3 \tau  \ge C N^\varepsilon$ for some $C >0 $ and $\varepsilon > 0$. 
Then 
\begin{align} &\notag %\label{:31a}& 
 \osc {N+l} ( \sqrt{2N} \cos \tau ) = 
\frac{1+ \Ord (N^{-1}) }{\sqrt{\pi \sin \tau }} 
\biggl(\frac{2}{N}\biggr)^{\frac{1}{4}} 
\\ & \notag 
\ts \biggl[ \sum_{n=0} ^{L-1} 
\sum_{m=0}^{n} C_{nm}^1 (N+l, \tau ) 
\sin \biggl\{ \frac{N}{2} (2 \tau - \sin 2\tau ) + 
D_{nm}^1 (\tau ) -(1+l)\tau \biggr\} +
\Ord (\frac{1}{N \sin \tau }) 
\biggr]
.\end{align}
\thetag{2} Let  $\tau >0 $. 
Assume that $N\sinh ^3 \tau  \ge C N^\varepsilon$ for some $C >0 $ and $\varepsilon > 0$. 
Then 
\begin{align}& \notag %\label{:31b}
 \osc {N+l} ( \sqrt{2N} \cosh \tau ) = 
\frac{1+ \Ord (N^{-1}) }{\sqrt{2\pi \sinh \tau }} 
\biggl(\frac{1}{2N}\biggr)^{\frac{1}{4}} 
\\ & \notag 
\ts \exp \biggl[\biggl(\frac{N+1+l}{2} \biggr)(2\tau -\sinh 2\tau ) + 
(1+l)\tau \biggr] 
\biggl[ \sum_{n=0} ^{L-1} 
\sum_{m=0}^{n} C_{nm}^2 (\tau , N+l) +
\Ord \biggl(\frac{\cosh ^3\tau }{N \sinh \tau } \biggr) \biggr]
.\end{align}
\end{lemma}
\PROOF 
\thetag{1} and \thetag{2} follow from \thetag{5.5} and \thetag{5.10} 
in \cite{KT11}, respectively. 
\PROOFEND

We next quote estimates from \cite{KT11,o-t.airy}. 
\begin{lemma}[\cite{KT11}, \cite{o-t.airy}] \label{l:32}
\thetag{1}  
Let $ y=\sqrt{2N} \cos \tau $ with $N \in \mathbb{N} $ and $0 < \tau \le \frac{\pi}{2} $. 
Assume that $N\sin ^3 \tau  \ge C N^\varepsilon$ for some $C >0 $ and $\varepsilon > 0$.
Then,
\begin{align}\notag %\label{:32a}
& \sum_{k=0}^{N-1} \osc {k}(y)^2 = 
 \frac{1}{\pi }\sqrt{2N-y^2} +\Ord \biggl(\frac{\rN }{2N-y^2 } \biggr)
.\end{align}
\thetag{2} 
Let $ y=\sqrt{2N} \cosh \tau $  with $N \in \mathbb{N} $ and $\tau > 0 $. 
Assume that $N\sinh ^3 \tau  \ge C N^\varepsilon$ for some $C >0 $ and $\varepsilon > 0$.
Then
\begin{align}\label{:32b}
 \sum_{k=0}^{N-1} \osc k (y)^2 =  \Ord \biggl(\frac{\rN }{y^2-2N}\biggr)
.\end{align}
\thetag{3} There is a positive constant $\Ct{;32a}$ such that for all $ N \in \N $
\begin{align}\label{:32c}&
\sup_{ y \in\R }| \osc N (y) | \le 
{\cref{;32a}}{N^{-\frac{1}{12}} } 
.\end{align}
\end{lemma}
\PROOF  
\thetag{1} follows from Lemma 5.2 \thetag{i} in \cite{KT11}. 
\thetag{2} follows from Lemma 5.2 \thetag{ii} in \cite{KT11}. 
From Lemma 6.9 in \cite{o-t.airy} there exists a constant $ \cref{;32a}$ such that 
\begin{align}& \notag %\label{:32d }
| N^{\frac{1}{12}} \osc N (2 \rN  + y N^{-\frac{1}{6}}) | \leq 
\frac{\cref{;32a}}{( 1 \lor | y | )^{\frac{1}{4}}} , 
\quad y \in [-2 N^{\frac{2}{3}} , \infty)
.\end{align}
Hence we have 
\begin{align}&\label{:32e}
| \osc N (y) | \le 
\frac{\cref{;32a}}
{
N^{\frac{1}{12}} 
( 1 \lor \{N^{\frac{1}{6}} |y - 2\rN  | \} )^{\frac{1}{4}}
}
 ,\quad y \in [0 , \infty)
.\end{align}
Claim \eqref{:32c} is immediate from \eqref{:32e} and the well-known property such that 
$ \osc N (y) = \osc N (-y) $ if $ N$ is even and that 
$ \osc N (y) = -\osc N (-y) $ if $ N$ is odd. 
\PROOFEND

\subsection{Determinantal kernels of $ \nN $-particle systems}
\label{s:32}

We recall the definition of determinantal point processes. 
Let $\map{K}{\sS  ^2}{\C}$ be a measurable kernel. 
A probability measure $ \mu $ on $ \SSS $ is called a determinantal point process 
with kernel $ K $ if, for each $ n $, 
its $ n $-point correlation function is given by 
\begin{align}& \label{:33l}
\rho ^n(x_1,\ldots,x_n) = \mathrm{det}[K(x_i,x_j)]_{i,j=1}^{n}
.\end{align}
If $ K $ is an Hermitian symmetric and of locally trace class such that 
$ 0 \le \mathrm{Spec} (K) \le 1 $, then there exists a unique determinantal point process with kernel $ K $ \cite{ST03,sos.00}. 

The distribution of the delabeled eigenvalues of GUE associated with \eqref{:10b} 
is a determinantal point process with kernel $\KN $ such that 
\begin{align}\label{:33m}
\KN (x,y) =\sum_{k=0}^{N-1} \osc{k}(x) \osc{k}(y)
.\end{align}
The Christoffel-Darboux formula and a simple calculation yield the following. 
\begin{align}&\label{:33n}
\KN (x,y)
=\sqrt{\frac{N}{2}} \frac{\osc N(x)\osc {N-1}(y)-\osc {N-1}(x)\osc N(y)}{x-y}
.\end{align}

From the scaling \eqref{:10e},  $ \muNt $ is a determinantal point process with kernel 
\begin{align}\label{:33o}&
\KtN (x,y) = \frac{1}{\sqrt{\nN }} 
\KN (\frac{x + \nN \theta }{\sqrt{\nN }}, \frac{y + \nN \theta }{\sqrt{\nN }}) 
.\end{align}
Let $ \xN = \rN x$ and $ \yN = \rN y $. We set 
\begin{align}\label{:33p}&
 \LNxy = \frac{1}{\sqrt{N}}\KN (\xN , \yN ) = 
 \frac{1}{\sqrt{N}} \KN (\sqrt{\nN } x , \sqrt{\nN } y ) 
.\end{align}
From \eqref{:33o} and \eqref{:33p} we then clearly see that 
\begin{align}\label{:33q}&
 \KtN (x,y) = \mathsf{L}^N  (\frac{x}{N} + \theta , \frac{y}{N} + \theta ) , 
\\ \notag &
 \LNxy = \KtN (N(x-\theta ), N (y- \theta )) 
.\end{align}
From \eqref{:33n} we deduce 
\begin{align} \label{:33r}& 
\LNx = \, (1/ \sqrt{2})
\{\osc {N-1}(\xN ) \osc {N}'(\xN ) - \osc {N}(\xN )\osc {N-1}'(\xN ) \} 
.\end{align}
Using the Schwartz inequality to \eqref{:33m} we see from \eqref{:33n} 
and \eqref{:33p} that 
\begin{align} \label{:33t}&
\LNyz ^2\le \LNy \LNz 
.\end{align}
From here on, we assume
\begin{align}\label{:33u}&
\aD < \alpha < \aU 
.\end{align}
We set 
\begin{align}\label{:33v}&
\BN = (-\sqrt{2}-N^{\alpha },-\sqrt{2}+N^{\alpha })
\cup (\sqrt{2}-N^{\alpha },\sqrt{2}+N^{\alpha })
.\end{align}
The next lemma will be used in \sref{s:4}. 
\begin{lemma}\label{l:33}
We set $ \UN = \R \backslash \BN $. Then the following holds. 
\\
\thetag{1} There exists a constant $ \Ct{;33a}$ such that for all $ N \in \N $ 
\begin{align}\label{:33x}&
\sup_{x,y \in \R }|\LNxy |  \le \cref{;33a}N^{\frac{1}{3}} 
,\\\label{:33z}& 
%\sup_{N\in \N }
\sup_{x,y \in \UN } |\LNxy | \le \cref{;33a} 
.\end{align}
\thetag{2} 
Assume \eqref{:33u}. Then there exists a constant $ \Ct{;33c}$ such that 
\begin{align}\label{:33a} 
 |\LNxy |\le  &\frac{\cref{;33c}}{{ N  |x-y| }}
\quad \text{ for each $ x,y\in \UN $, $ N \in \N $}
.\end{align}
\end{lemma}
\PROOF  It is well known that 
\begin{align}\notag %\label{:33d}
\sqrt{2} \osc {n}'(x)=\sqrt{n} \osc {n-1} (x) -\sqrt{n+1}\osc {n+1}(x)
.\end{align}
From this and \eqref{:33r}, we see that with a simple calculation  
\begin{align}&\label{:33e}
\LNx =
\frac{1}{\sqrt{2}} \{ \osc {N-1} \osc {N}' - \osc {N}\osc {N-1}' \} (\xN )
\\\notag &  = \frac{N^{\frac{1}{2}}}{2}\{ 
\osc {N-1}^2 + \osc {N}^2 - 
{\sqrt{1 - N^{-1}}} \osc {N-2}\osc {N} - 
{\sqrt{1 + N^{-1}}} \osc {N-1}\osc {N+1} \} (\xN ) 
.\end{align}
Combining this with \eqref{:32c} we obtain 
\begin{align}& \notag 
 \LNx \le 
\frac{N^{\frac{1}{2}} }{2} 5\cref{;32a}^2 N^{-\frac{1}{6}} =
\frac{5\cref{;32a}^2}{2}  N^{\frac{1}{3}}
.\end{align}
From this and \eqref{:33t} we deduce \eqref{:33x}.  
From \lref{l:31} and \eqref{:33e}, we see that 
\begin{align}&\notag 
\sup_{N\in \N }\sup_{y\in \UN } \LNy < \infty 
.\end{align}
We deduce \eqref{:33z} from this and \eqref{:33t}. 
Taking a constant $ \cref{;33a}$ in \eqref{:33x} and \eqref{:33z} in common 
 completes the proof of \thetag{1}. 

Claim \eqref{:33a} follows from \lref{l:31}, \eqref{:33n}, and \eqref{:33p}. 
\PROOFEND

\section{Proof of \eqref{:26a}--\eqref{:26d}}\label{s:4}

As we see in \sref{s:2}, the point of the proof of \tref{l:11} is to check 
conditions \eqref{:26a}--\eqref{:26d} in \lref{l:26}. 
The purpose of this section is to prove these equations. 
We recall a property of the reduced Palm measures of determinantal point processes. 
\begin{lemma}[\cite{ST03}] \label{l:41}
Let $\mu $ be a determinantal point process with kernel $K$. 
Assume that $ K(x,y) = \overline{K(y,x)} $ and $ 0 \le \mathrm{Spec} (K) \le 1 $. 
Then the reduced Palm measure $\mu _x$ is 
a determinantal point process with kernel $K_x$ given by 
\begin{align}\label{:41x}
K_x(y,z)= K(y,z) - \frac{K(y,x)K(x,z)}{K(x,x)}
\end{align}
for $ x $ such that $ K(x,x)>0 $. 
\end{lemma}

Let $ \KtN $  be the determinantal kernel of $ \muNt $ given by \eqref{:33o}. 
Let $  \muNtx$ be as in \lref{l:26}. 
Recall that $ \KtN (y,z) = \KtN (z,y) $ by definition. 
Then from this, \eqref{:33o}, and \eqref{:41x}, 
 $  \muNtx $ is a determinantal point process with kernel 
\begin{align}\label{:41y}&
\mathsf{K}_{\theta ,x}^{\nN } (y,z) = 
\KtN (y,z) - \frac{ \KtN (x,y) \KtN (x,z) }{\KtN (x,x)}
.\end{align}
From \eqref{:33l} and \eqref{:41y}, we calculate correlation functions 
 in \eqref{:26a}--\eqref{:26d} as follows. 
\begin{align}\label{:41a}&
\rho_\theta ^{N,1}(y) = \KtN (y,y)
,\\\label{:41b}& 
\rho_{\theta ,x}^{N,1}(y)  -\rho_\theta ^{N,1}(y)  = - 
\frac{\KtN (x,y) ^2}{\KtN (x,x)} 
,\\\label{:41c}& 
\rho_{\theta }^{N,2}(y,z)-\rho_\theta ^{N,1}(y) \rho_\theta ^{N,1} (z) 
= 
 -\KtN (y,z)^2 
,\\\label{:41d}&
\rho_{\theta ,x}^{N,2}(y,z)-\rho_{\theta ,x} ^{N,1}(y) \rho_{\theta ,x} ^{N,1} (z)
- \{\rho_\theta ^{N,2}(y,z) -\rho_\theta ^{N,1}(y) \rho_\theta ^{N,1} (z) \} 
\\ \notag & \quad \quad 
 = 
 - \mathsf{K}_{\theta ,x}^{\nN }(y,z)^2 + \KtN (y,z)^2 
\\ \notag & \quad \quad 
=
 2\frac{\KtN (y,z) \KtN (x,y) \KtN (x,z)  }{\KtN (x,x)} 
 - 
\frac{  \KtN (x,y)^2 \KtN (x,z) ^2}{\KtN (x,x)^2 }  
.\end{align}
Using these and \eqref{:33q} we rewrite \eqref{:26a}--\eqref{:26d} as follows. 
\begin{lemma} \label{l:42}
To simplify the notation, let 
\begin{align} \label{:42z} &
 \xxN = \frac{x}{\nN } + \theta , \quad 
\Trx = \{ y \in \R \, ;\, \frac{r}{N} \le | \xxN - y | < \infty \} 
.\end{align}
Then \eqref{:26a}--\eqref{:26d} are equivalent to 
\eqref{:42A}--\eqref{:42D} below, respectively. 
\begin{align} 
& \label{:42A}
\1 \Sup{ \int_{\Trx } \frac{\LNy }{\ellxy } dy - \theta   }   =0
,\\& \label{:42B}
\1 \Sup{ \int_{\Trx } \5  \frac{\LNyxN ^2}{\LNxN } dy  }  =0
.\end{align}
%Furthermore, 
\begin{align} & \label{:42C}
\1 \Sup{ \int_{\Trx } \frac{1}{\nN } \frac{\LNy }{\ellxyA ^2} dy
-\int_{\Trxx } \frac{\LNyz ^2 }{\ellxyP \ellxzP } dydz  }  =0
,\\
&\label{:42D}
\1 \Sup{\int_{\Trx } \frac{1}{\nN }\frac{1}{\ellxyA ^2}  
\frac{\LNyxN ^2}{\LNxN }  dy 
\\& 
\notag 
+\int_{\Trxx } \frac{1 }{\ellxyP \ellxzP }
\\ \notag & \quad \quad \quad 
\biggl\{ 2\frac{\0 }{\LNxN } - 
\frac{ \LNxyN \LNxz }{\LNxN ^2 }  \biggr\}  dydz }  =0 
.\end{align}
\end{lemma}
\PROOF 
Recall that 
$ \LNxy = \KtN (N(x-\theta ), N (y- \theta )) $ by \eqref{:33q}.  
 Then \lref{l:42} follows easily from \eqref{:41a}--\eqref{:41d}. 
\PROOFEND 

Let $ \BN $ and $ \UN  $ be as in \lref{l:33}. 
Decompose $ \UN   $ into $ \UNone  $ and $ \UNtwo $  such that 
\begin{align} \notag %\label{:43z} 
 \UNone  &= [ - \sqrt{2}+N^{\alpha }, \sqrt{2}-N^{\alpha } ] 
,\quad  %\\ \notag 
 \UNtwo  = \R \backslash  ( - \sqrt{2}-N^{\alpha }, \sqrt{2}+N^{\alpha } )
.\end{align}
Then clearly $ \UN  = \UNone  \cup \UNtwo  $ and 
$ \R  = \UNone \cup \UNtwo  \cup \BN $. 
We begin by the integral outside $ \UNone  $. 
\begin{lemma} \label{l:43} 
Let $ 0 < q < 3/2 $. Then 
\begin{align}& \label{:43a}
\limsupi{N} 
\Sup{  \int_{\R \backslash \UNone } \frac{ \LNy ^q}{ | \ellxy | }   dy } = 0 
.\end{align}
\end{lemma}
\PROOF  
From \eqref{:33x}, \eqref{:42z}, and the definition of $ \BN $, we obtain that 
\begin{align} \label{:43c}
&\limsupi{N} 
 \Sup{  \int_{ \BN  } \frac{ \LNy ^q }{|\ellxy |}  dy  }  
 \\\notag
\le & \limsupi{N} 
 \Sup{ \int_{ \BN  } \frac{\cref{;33a} ^q  N^{\frac{q}{3}} }{| \ellxy |} dy } 
\\\notag 
\le &  \limsupi{\nN }
\Sup{\cref{;33a} ^q    
N^{\frac{q}{3}} 
\Bigl\{ \log \Bigl| \frac{x}{N} +\theta - (\sqrt{2} - N^{\alpha }) \Bigr| 
- \log \Bigl| \frac{x}{N} +\theta - (\sqrt{2} + N^{\alpha} ) \Bigr| \Bigr\} 
\quad \quad \quad {}
\\ \notag  \quad \quad \quad \quad 
+ &\cref{;33a}  ^q 
N^{\frac{q}{3}} 
\Bigl\{ \log \Bigl| \frac{x}{N} +\theta - (-\sqrt{2} - N^{\alpha }) \Bigr| 
- \log \Bigl| \frac{x}{N} +\theta - (-\sqrt{2} + N^{\alpha} ) \Bigr| \Bigr\} 
}
\\ \notag 
=& \mathcal{O}(N^{\frac{q}{3}+ \alpha }) = 0 \quad \text{ as $ N \to \infty $}
.\end{align}
Here we used $ q < 3/2 $ and $  \alpha < - 1/2 $ in the last line. 

Note that 
$ | y | \ge \sqrt{2} + N^{\alpha } $ for $y \in \UNtwo $. 
Let $ y=\sqrt{2} \cosh \tau $. Then we see that 
\begin{align}\notag %\label{:43d}
N \sinh ^3 \tau &= N(\cosh ^2 \tau -1)^{\frac{3}{2}} 
\\& \notag
=N2^{-\frac{3}{2}}(y^2-2)^{\frac{3}{2}} 
%\\& \notag
\ge N2^{-\frac{3}{2}}(2\sqrt{2}N^{\alpha }+N^{2\alpha })^{\frac{3}{2}} 
.\end{align}
From this, $ q > 0 $, and $\alpha >-{2}/{3}$, %and \eqref{:43d}, 
 we apply \eqref{:32b} to obtain $\Ct{;19}> 0 $ such that,
\begin{align} &\notag %\label{:43e} 
\limsupi{N} \Sup{  \int_{ \UNtwo  } \frac{\LNy ^q }{|\ellxy |}  dy  }  
\le  \limsupi{N} \Sup{ \int_{ \UNtwo  } 
 \frac{\cref{;19}}{ |\ellxy | N ^q (y^2-2) ^q } dy } = 0
,\end{align}
which combined with \eqref{:43c} yields \eqref{:43a}. 
\PROOFEND

\begin{lemma} \label{l:44}
\eqref{:42A} holds. 
\end{lemma}

\PROOF  Let $ y=\sqrt{2} \cos \tau $. Then 
$ N \sin ^3 \tau \ge N 2^{-\frac{3}{2}} (2\sqrt{2}N^{\alpha }-N^{2\alpha })$ 
for $y\in \UNone $.  Then applying \lref{l:32} (1) we deduce that for each $ r > 0 $
\begin{align}&\notag %\label{:44a}
\limsupi{N} 
\Sup{ \int_{ \Trx \cap \UNone  } \frac{\LNy }{\ellxy } dy - \theta  }   
\\\notag
= &\limsupi{N} 
 \Sup{ \Big\{ \int_{- \sqrt{2} + N^{\alpha } } ^{\xxN - \frac{r}{N} }
 + 
\int_{\xxN + \frac{r}{N}}^{\sqrt{2} - N^{\alpha } } \Big\} 
\frac{1}{\ellxy }   
\frac{1}{\pi }\sqrt{2-y^2} \,  dy - \theta   }  
\\\notag
=&\biggl| \mathrm{P.V.} \int_{-\sqrt{2} }^{\sqrt{2}} 
 \frac{1}{\theta -y } \frac{1}{\pi }\sqrt{2-y^2}  dy -\theta  \biggr|=0
.\end{align}
Combining this with \eqref{:43a}, we obtain \eqref{:42A}. 
\PROOFEND 

It is well known that $ \KtN (x,x)$ are positive and continuous in $ x $, and 
$ \{ \KtN (x,x) \}_{N\in\N } $ converges to 
$ \mathsf{K}_{\theta } (x,x)=\sqrt{2-\theta ^2}/\pi  $ 
uniformly on each compact set. Then we see 
\begin{align}& \notag %\label{:45y}&
\sup_{N \in \N }\sup_{x \in S_R } \frac{1 }{\KtN (x,x)} < \infty 
.\end{align}
From this, \eqref{:33q}, and \eqref{:42z}, we see that 
the following constant $ \Ct{;44b}$ is finite. 
\begin{align}\label{:45z}&
\cref{;44b}:= \sup_{N \in \N} \sup_{x\in S_R }\frac{1}{\LNxN } < \infty 
.\end{align}
\begin{lemma} \label{l:45} \eqref{:45a} and \eqref{:45A} below hold. 
In particular, \eqref{:42B} holds.
\begin{align}\label{:45a}&
\1 \Sup{ \int_{\Trx } \frac{\LNxyN ^2}{\ellxyA \LNxN } 
dy  }  =0
,\\ &\label{:45A}
\1
\Sup{ \int_{\Trx } 
 \frac{\LNxyN }{\ellxyA \LNxN }  dy } = 0 
.\end{align}
\end{lemma}
\PROOF  
From \eqref{:33t} and \eqref{:43a} we deduce that as $ N \to \infty $
\begin{align}\label{:45b}&
 \Sup{  \int_{\R  \backslash \UNone   }  \frac{\LNxyN ^2}{\ellxyA \LNxN }  dy }
%\\ \notag &
\le  \Sup{  \int_{\R  \backslash \UNone   }  \frac{\LNy }{\ellxyA } dy } \to 0 
.\end{align}
From \eqref{:33a} and \eqref{:45z} for each $ N \in \N $ and $ r > 0 $
\begin{align}\label{:45c}
\SSup{\int_{ \Trx \cap \UNone }  \frac{\LNxyN ^2 \, dy }{\ellxyA \LNxN } } 
 & \le 
\SSup{ \int_{ \Trx \cap \UNone } 
\frac{\cref{;33c} ^2 \cref{;44b} \, dy }{N^2 |\xxN -y|^3}}
\\\notag & \le \frac{\cref{;33c} ^2 \cref{;44b}}{r^2} 
.\end{align}
Hence \eqref{:45a} follows from \eqref{:45b} and \eqref{:45c}. 
This completes the proof of \eqref{:45a}. 

We next prove \eqref{:45A}.  
From \eqref{:33t}, \eqref{:43a}, and \eqref{:45z} we see for each $ r > 0 $ 
\begin{align}\label{:45e}&
\limsupi{N} 
\Sup{\int_{\Trx \backslash \UNone } \frac{\LNxyN }{\ellxyA \LNxN } 
  dy }  =0 
.\end{align}
From \eqref{:33a} and \eqref{:45z} 
we see that for each $ N \in \N $ and $ r > 0 $
\begin{align}\label{:45f} 
 \Sup{ \int_{\Trx \cap \UNone } 
\frac{\LNxyN  \, dy }{\ellxyA \LNxN } } 
\le &
\Sup{\int_{\Trx \cap \UNone } 
\frac{\cref{;33c}\cref{;44b}\,  dy  }{N \ellxyA ^2} } 
\\\notag \le  & 
\frac{2\cref{;33c}\cref{;44b}}{r}
.\end{align}
Combining \eqref{:45e} and \eqref{:45f} we obtain \eqref{:45A}.  
\PROOFEND 

\begin{lemma}\label{l:46}
\eqref{:46a} and \eqref{:46b} below hold. In particular, \eqref{:42C} holds. 
\begin{align}\label{:46a}&
\1 \Sup{ \int_{\Trx }  \frac{\LNy }{ \nN \ellxyA ^2} \,dy  }  =0
,\\\label{:46b}&
\1 \Sup{ \int_{\Trxx }\frac{\LNyz ^2 }{\ellXYZ } dydz }  
=0
.\end{align}
\end{lemma}
\PROOF 
Note that $ \LNy \le \cref{;33a} $ on $ \UN $ by \eqref{:33z}.  
Then by the definition of $ \Trx $, 
\begin{align}\label{:46c}&
 \int_{\Trx \cap \UN  }  \frac{\LNy }{ \nN \ellxyA ^2} \,dy  \le 
\frac{\cref{;33a} }{N} \frac{2N}{r} = \frac{2\cref{;33a}}{r}
.\end{align}
By \eqref{:33x} we see $ \LNy \le \cref{;33a}N^{\frac{1}{3}}$ on $ \R $. 
Recall that $ |\BN | = 4 N^{\alpha }$ by construction. 
Furthermore, 
$ \Ct{;46d}:=\limsupi{N} \sup_{y \in \BN }\| {\ellxyA ^{-2}} \|_{R} < \infty $. 
Hence for each $ r > 0 $ 
\begin{align}\label{:46d}&
\limsupi{N}\int_{\Trx \cap \BN  }  \frac{\LNy }{ \nN \ellxyA ^2} \,dy 
 \le \limsupi{N} 
\frac{ \cref{;33a} N^{\frac{1}{3}} 4 N^{\alpha }\cref{;46d}}{N} 
 = 0 
.\end{align}
Here we used $ \alpha < -1/2$. 
We thus obtain \eqref{:46a} from \eqref{:46c} and \eqref{:46d}. 

We proceed with the proof of \eqref{:46b}.  
We first consider the integral away from the diagonal line. 
By \eqref{:33a} and the Schwartz inequality, we see that 
\begin{align}\notag %\label{ :46e}
& \quad \ 
\Sup{ \int_{(\Trx \cap \UN )^2 \cap \{ |y-z| \ge \qN \} } 
\frac{ \LNyz ^2 }{\ellXYZ } dydz } 
\\ \notag &
\le \Sup{ 
\int_{(\Trx \cap \UN )^2 \cap \{ |y-z| \ge \qN \} } 
\frac{ \cref{;33c} ^2 } {N^2 |y-z|^2 |\xxN -y | |\xxN -z |  }dydz } 
\\ \notag &
\le \Sup{ 
\Big\{ \int_{\Trxx  \cap \{ |y-z| \ge \qN \} } 
 \frac{ \cref{;33c} ^2 } {N^2 |y-z|^2 |\xxN -y |^2  }dydz \Big\}^{\frac{1}{2}}
\\ \notag & \quad \quad \quad \quad \quad 
\Big\{ \int_{\Trxx  \cap \{ |y-z| \ge \qN \} } 
 \frac{ \cref{;33c} ^2 } {N^2 |y-z|^2 |\xxN -z |^2  }dydz \Big\}^{\frac{1}{2}}
 } 
\\ \notag &
= \Sup{ \int_{\Trxx  \cap \{ |y-z| \ge \qN \} } 
 \frac{ \cref{;33c} ^2 } {N^2 |y-z|^2 |\xxN -y |^2  }dydz }
\\ \notag &
\le \cref{;33c} ^2\frac{2N}{N^2}\Big\{\frac{2N}{r}\Big\} = 
\frac{4\cref{;33c} ^2}{r} 
.\end{align}
The last line follows from a straightforward calculation. 
Indeed, first integrating $ z $ over $ \{ |y-z| \ge \frac{1}{N}\}  $, and then integrating $ y$ over $ \Trx $, we obtain the inequality in the last line. We therefore see that 
\begin{align}\label{:46j}&
\limi{r} \limi{N} \Sup{ \int_{(\Trx \cap \UN )^2 \cap \{ |y-z| \ge \qN \} } 
\frac{ \LNyz ^2 }{\ellXYZ } dydz } 
 = 0 
.\end{align}

We next consider the integral near the diagonal. 
From \eqref{:33z}, we see that 
\begin{align}\label{:46k}&
\Sup{ \int_{(\Trx \cap \UN )^2 \cap \{ |y-z| \le \qN \} } 
\frac{ \LNyz ^2 }{\ellXYZ } dydz } 
\\ \notag \le & 
\Sup{ \int_{(\Trx \cap \UN )^2 \cap \{ |y-z| \le \qN \} } 
\frac{ \cref{;33a}^2 }{\ellXYZ } dydz } 
\\ \notag \le & 
\Sup{ \int_{\Trxx  \cap \{ |y-z| \le \qN \} } 
\frac{\cref{;33a}^2}{2} \{
 \frac{ 1}{\ellxyA ^2 }+  \frac{ 1}{\ellxzA ^2 }
\} 
 dydz }
\\ \notag 
= & 
\frac{2\cref{;33a}^2}{ N } \Sup{\int_{\Trx } \frac{1}{\ellxyA ^2 }dy } 
= 
 \frac{2\cref{;33a}^2}{ N } \frac{2N}{r} = \, \frac{4\cref{;33a}^2}{r}
.\end{align}
From \eqref{:46j} and \eqref{:46k}, we have 
\begin{align}\label{:46l}&
\limi{r} \limi{N} 
 \Sup{ \int_{(\Trx \cap \UN )^2} 
\frac{ \LNyz ^2 }{\ellXYZ } dydz } 
= 0 
.\end{align}

We next consider the integral on $ \BN \ts \BN $. Let 
$$ \Ct{;45b}=\limsupi{N}\sup_{x\in \SR , y \in \BN }|\ellxy |^{-1} .$$
Then, we deduce from \eqref{:33x} and the definition of $ \BN $ 
given by \eqref{:33v} that 
\begin{align}\label{:46m}
\limsupi{N} &\Sup{ \int_{(\Trx \cap \BN )^2} 
\frac{ \LNyz ^2 }{\ellXYZ } dydz } \\ \notag 
\le &\limi{N} \cref{;33a}^2 \cref{;45b}^2 N^{\frac{2}{3}} (4N^{\alpha })^2 
= \,  0 
.\end{align}
Here we used $ |\BN |= 4N^{\alpha }$ for the inequality and $ \alpha < - 1/2 $ for the last equality. 

We finally consider the case $ \UN \ts \BN $. 
Then a similar argument yields 
\begin{align}\label{:46n}&
%\limsupi{N} 
\Sup{ \int_{(\Trx \cap \UN )\ts (\Trx \cap \BN )} 
\frac{ \LNyz ^2 }{\ellXYZ } dydz } 
\\ \notag 
\le & % \limsupi{N} 
\Sup{ \int_{(\Trx \cap \UN )\ts (\Trx \cap \BN )} 
\frac{ \LNy \LNz }{\ellXYZ } dydz } 
\\ \notag 
=& % \limsupi{N} 
\Sup{ \int_{\Trx \cap \UN } 
\frac{ \LNy }{\ellxyA } dy 
\int_{ \Trx \cap \BN } 
\frac{ \LNz }{\ellxzA } dz } 
\\ \notag 
=& \mathcal{O}(\log N ) \mathcal{O}(N^{\frac{1}{3} + \alpha })  
\to 0 \quad \text{ as } N \to \infty 
.\end{align}

Collecting \eqref{:46l}, \eqref{:46m}, and \eqref{:46n}, we conclude \eqref{:46b}.  
\PROOFEND 

\begin{lemma}\label{l:47}
 \eqref{:42D} holds.
\end{lemma}
\PROOF 
We shall estimate the three terms in \eqref{:42D} beginning with the first. 
From \eqref{:33t} and \eqref{:46a} we have %as $ N \to \infty $
\begin{align}\label{:47a}&
\1 
\Sup{\int_{\Trx } 
\frac{\LNyxN ^2 dy }{\nN \ellxyA ^2 \LNxN } }
\\ \notag 
\le &
\1 
\Sup{\int_{\Trx } \frac{\LNy  dy }{\nN \ellxyA ^2} } 
=  0 
.\end{align}
Next, using the Schwartz inequality, we have for the second term 
\begin{align}&\notag % \label{:47g}&
%\1 
\Sup{\int_{\Trxx } 
\frac{ \0  \, dydz }{\ellxyA \ellxzA \LNxN }}
\\ \notag \le &
\Sup{
\int_{\Trxx } 
 \frac{ \LNyz ^2 dydz  }{ \ellxyA \ellxzA  } } ^{\frac{1}{2}} 
\Sup{\int_{\Trxx } 
 \frac{ \LNxyN ^2 \LNxz ^2  dydz }{ \ellxyA \ellxzA  \LNxN ^2 }}^{\frac{1}{2}}
%\\ \notag &
\\ \notag = &
\Sup{\int_{\Trxx  } 
 \frac{ \LNyz ^2  dydz }{ \ellxyA \ellxzA } }^{\frac{1}{2}} 
\Sup{\int_{\Trx }
 \frac{ \LNxyN ^2 }{|\ellxy | \LNxN  } dy  
}
.\end{align}
Applying \eqref{:46b} and \eqref{:45a} to the last line, we obtain 
\begin{align}\label{:47b}&
\1 
\Sup{\int_{\Trxx } 
\frac{ \0  \, dydz }{\ellxyA \ellxzA \LNxN }}= 0 
.\end{align}
We finally estimate the third term. From \eqref{:45A}, as $ N \to \infty $, we have 
\begin{align}\label{:47c}&
\Sup{\int_{\Trxx } 
 \frac{ \LNxyN \LNxz \, dydz }{\ellxyA \ellxzA \LNxN ^2 }}
\\ \notag  
= & 
\Sup{ \Big\{ 
\int_{\Trx } \frac{\LNxyN \, dy }{\ellxyA \LNxN } \Big\}^2  }  
\\ \notag   = &%\1
\Sup{ \int_{\Trx } \frac{\LNxyN \,   dy }{\ellxyA \LNxN } }^2  \to 0 
\quad \text{ by \eqref{:45A}}
.\end{align}
From \eqref{:47a}, \eqref{:47b}, and \eqref{:47c} we obtain \eqref{:42D}.  
This completes the proof.     
\PROOFEND

\section{Proof of \tref{l:11}} \label{s:5} 
From \lref{l:44}--\lref{l:47} we deduce that all the assumptions 
\eqref{:26a}--\eqref{:26d} in \lref{l:26} are satisfied. Hence 
\eqref{:23k} is proved by \lref{l:26}. 
Then \tref{l:11} follows from \lref{l:24}, \lref{l:25}, and \lref{l:26}.

\section{Proof of \tref{l:12}}  \label{s:6}
In this section we prove \tref{l:12} using \tref{l:11}. 
It is sufficient for the proof of \tref{l:12} to prove \eqref{:11f} in $ C([0,T];\R ^m)$ for each $ T \in \N $. Hence we fix $  T \in \N $.  
Let $ \mathbf{X}^{\nN }=(X^{N , i})_{i=1}^{\nN }$ be as in \eqref{:11c}.  
%a solution of \eqref{:11c} as before. 
Let $   \Y ^{\theta , N,i}  = \{ \Y _t^{\theta , N,i} \} $ such that 
\begin{align}\label{:62p}&
\Y _t^{\theta , N,i} = X_t^{N , i} + \theta t 
.\end{align}
Then from \eqref{:11c} we see that 
$ \mathbf{\Y }^{\theta ,N} = (\Y ^{\theta ,N,i})_{i=1}^N$ 
is a solution of 
\begin{align}\label{:62q} \quad 
d \Y _t^{\theta ,N,i} = \, &dB_t^{i} + \sum_{ j \neq i }^{\nN }
 \frac{1}{\Y _t^{\theta ,N,i}  - \Y _t^{\theta ,N,j} }dt 
 - \frac{1}{N}\Y _t^{\theta ,N,i}\,dt 
 + 
\frac{\theta }{N} 
\,dt 
\end{align}
with the same initial condition as $ \mathbf{X}^{N}$. 
Let $ P^{\theta , N}$ and $ Q^{\theta , N}$ be 
the distributions of $ \mathbf{X}^{N}$ and $ \mathbf{\Y }^{\theta ,N}$ on 
$ C([0,T];\R ^N )$, 
respectively. Then applying the Girsanov theorem \cite[pp.190-195]{IW} 
to \eqref{:62q}, we see that 
\begin{align}\label{:62r}
\frac{d Q^{\theta , N}}{d P^{\theta , N}} (\mathbf{W}) = &
\exp \{ \int_0^T \sum_{i=1}^N \frac{\theta }{N} dB_t^i 
-\frac{1}{2}
\int_0^T \sum_{i=1}^N \Big| \frac{\theta }{N}\Big|^2 dt 
\} 
\\ \notag 
 = &
\exp \{  \frac{\theta }{N}\sum_{i=1}^N B_T^i 
- \frac{\theta ^2 T }{2N } \} 
,\end{align}
where we write $ \mathbf{W}=(W^i)\in C([0,T];\R ^N )$ and 
$ \{ B^i \}_{i=1}^N $ under $  P^{\theta , N}$ 
are independent copies of Brownian motions starting at the origin. 
\begin{lemma} \label{l:62}
For each $  \epsilon > 0 $, 
\begin{align}\label{:62z}&
\limi{N}
Q^{\theta , N} \Big( \Big | 
\frac{d P^{\theta , N}}{d Q^{\theta , N}} (\mathbf{W})-1 \Big | \ge \epsilon  \Big ) = 0 
.\end{align}
\end{lemma}

\PROOF 
It is sufficient for \eqref{:62z} to prove, for each $  \epsilon > 0 $, 
\begin{align}\notag %\label{:62b}
&
\limi{N}
P^{\theta , N} \Big( \Big | 
\frac{d Q^{\theta , N}}{d P^{\theta , N}} (\mathbf{W})-1 \Big | \ge \epsilon  \Big ) = 0 
.\end{align}
This follows from \eqref{:62r} immediately. 
\PROOFEND

\noindent {\em Proof of \tref{l:12}. } 
We write $ \mathbf{W}^m =(W^1,\ldots,W^m) \in C([0,T];\R ^m )$ for $ \mathbf{W}=(W^i)_{i=1}^{N}$, where $ m \le N \le \infty $. 
Let  
$ Q^{\theta }$ be the distribution of the solution $ \mathbf{\Y }^{\theta }$ with initial distribution $ \mut \circ \lab ^{-1}$. 
From \tref{l:11} and \eqref{:62p}  we deduce that for each $ m \in \N $ 
\begin{align}& \notag %\label{:62b}&
\limi{N}Q^{\theta , N}( \mathbf{W}^m \in \cdot ) = Q^{\theta }( \mathbf{W}^m \in \cdot )
\end{align}
weakly in $ C([0,T];\mathbb{R}^m)$. 
Then from this, for each 
$ F \in C_{b}(C([0,T];\mathbb{R}^m))$,  
\begin{align}\label{:64b}&
\limi{N} \int_{C([0,T];\mathbb{R}^N)} F (\mathbf{W}^m)dQ^{\theta , N}= 
\int_{C([0,T];\mathbb{R}^{\N })} F (\mathbf{W}^m)dQ^{\theta }
.\end{align}
We obtain from \eqref{:62z} and \eqref{:64b} that 
\begin{align}\notag %\label{:64c}
\limi{N}\int_{C([0,T];\mathbb{R}^N)} F (\mathbf{W}^m)dP^{N , \theta }
=&
\limi{N}\int_{C([0,T];\mathbb{R}^N)} F (\mathbf{W}^m)
\frac{d P^{\theta , N}}{d Q^{\theta , N}} (\mathbf{W})
dQ^{\theta ,N }
\\ \notag 
=&
\limi{N}\int_{C([0,T];\mathbb{R}^N)} F (\mathbf{W}^m) dQ^{\theta ,N }
\\ \notag 
 = &
\int_{C([0,T];\mathbb{R}^{\N })} F (\mathbf{W}^m)dQ^{\theta }
.\end{align}
This implies \eqref{:11f}.  
We have thus completed the proof of \tref{l:12}.  
\qed

%\section{Acknowledgments} \label{s:6}
\section{Acknowledgement}
H.O.\! thanks Professor H. Spohn for a useful comment at RIMS in Kyoto University in 2002. H.O.\! is supported in part by JSPS KAKENHI Grant Nos. 16K13764, 
16H02149, 16H06338, and KIBAN-A 24244010. 
Y.K. is supported by Grant-in-Aid for JSPS JSPS Research Fellowships (Grant No.\!\! 15J03091).
% \end{acknowledgement}
%\input{sg-ref.tex}

\noindent 
Yosuke Kawamoto \\
Faculty of Mathematics, Kyushu University, Fukuoka 819-0395, Japan
\\
\texttt{y-kawamoto@math.kyushu-u.ac.jp}

\bs
\noindent 
Hirofumi Osada \\
Faculty of Mathematics, Kyushu University, Fukuoka 819-0395, Japan
\\
\texttt{osada@math.kyushu-u.ac.jp}
\\


\begin{thebibliography}{99}
\DN\ttl[1]{{#1},}
% \DN\ttl[1]{\textit{#1},} 
%\DN\ttj[1]{\textit{#1},}  %%for references %%

\bibitem{AGZ}
Anderson, G.W., Guionnet, A., Zeitouni, O., \ttl{An Introduction to Random Matrices} Cambridge university press, 2010.

\bibitem{bey.14-duke} Bourgade, P., Erd\"{o}s, L., and Yau, H.-T., \ttl{Universality of general $ \beta $-ensembles} Duke Math. J. {\bf 163},  (2014), 1127-1190. 

%\bibitem{Dys62}Dyson, F. J. :\ttl{A Brownian-motion model for the eigenvalues of a random matrix} J. Math. Phys. {\bf 3}, 1191-1198 (1962).

%\bibitem{forrester} Forrester, Peter J., \ttl{Log-gases and Random Matrices} London Mathematical Society Monographs, Princeton University Press (2010).

%\bibitem{Fr}Fritz, J., \ttl{Gradient Dynamics of Infinite Point Systems} Ann. Probab. {\bf 15} (1987) 478-514. 

%\bibitem{fot} Fukushima, M., {\it et al.}, \ttl{Dirichlet forms and symmetric Markov processes} 2nd ed., Walter de Gruyter (2011). 

%\bibitem{o-h.bes} Honda,~R., Osada,~H, \ttl{Infinite-dimensional stochastic differential equations related to Bessel random point fields}  (Published on line: Stochastic Processes and Their Applications, http://dx.doi.org/10.1016/j.spa.2015.05.005) arXiv:1405.0523 [math.PR]. 


\bibitem{IW} Ikeda, N., Watanabe, S., \ttl{Stochastic differential equations and diffusion processes} 2nd ed, North-Holland (1989). 

\bibitem{inu} Inukai, K., \ttl{Collision or non-collision problem for interacting Brownian particles} Proc. Japan Acad. Ser. A {\bf 82}, (2006), 66-70.

%\bibitem{KS} Karatzas, I., Shreve, E., \ttl{Brownian motion and stochastic calculus} 2nd ed, Springer (1991)

%\bibitem{KT07b} Katori, M., Tanemura, H.: \ttl{Noncolliding Brownian motion and determinantal processes}. J. Stat. Phys. {\bf 129}, 1233-1277 (2007). 

%\bibitem{KT09} Katori, M., Tanemura, H.: Zeros of Airy function and relaxation process, {\it J. Stat. Phys.} {\bf 136}, 1177--1204 (2009). 

\bibitem{KT10a}
Katori, M., Tanemura, H., \ttl{Non-equilibrium dynamics of Dyson's model with an infinite number of particles} Comm. Math. Phys. {\bf 136}, 1177--1204 (2010)

\bibitem{KT11}Katori, M., Tanemura, H., \ttl{Markov property of determinantal processes with extended sine, Airy, and Bessel kernels} Markov processes and related fields {\bf 17}, 541-580 (2011). 

\bibitem{kawa.dp} Kawamoto, Y., \ttl{Density preservation of unlabeled diffusion in systems with infinitely many particles} in Stochastic analysis on large scale interacting systems, 337--350, RIMS K\^oky\^uroku Bessatsu, B59, Res. Inst. Math. Sci. (RIMS), Kyoto, 2016. 

\bibitem{k-o.fpa} Kawamoto, Y., Osada, H., \ttl{Finite-particle approximations for interacting Brownian particles with logarithmic potentials} to appear in J. Math. Soc. Japan, arXiv:1607.06922v1 [math.PR] 
% 
% \bibitem{k-o.sgap} Kawamoto, Y., Osada, H.: Universality of Dyson's model in infinite dimensions and SDE gaps. (in preparation)
% 
% \bibitem{k-o.sgap} Kawamoto, Y., Osada, H.: Dynamical bulk scaling of Gaussian unitary ensembles and SDE gaps. (in preparation) 第2案


%\bibitem{lang.1} Lang,~R., \ttl{Unendlich-dimensionale Wienerprocesse mit Wechselwirkung I}  Z. Wahrschverw. Gebiete  {\bf  38  } (1977) 55-72.   

%\bibitem{lang.2} Lang,~R., \ttl{Unendlich-dimensionale Wienerprocesse mit Wechselwirkung II} Z. Wahrschverw. Gebiete  {\bf  39  } (1978) 277-299.   

%\bibitem{mr} Ma, Z.-M. and R\"ockner, M., \ttl{Introduction to the theory of (non-symmetric) Dirichlet forms} \ttl{Springer-Verlag} 1992. 

\bibitem{mehta} Mehta, M. L., \ttl{Random Matrices} 3rd edition, Amsterdam: Elsevier, 2004. 

%\bibitem{o.dfa} Osada, H., \ttl{Dirichlet form approach to infinite-dimensional Wiener processes with singular interactions} Commun. Math. Phys. {\bf 176}, 117-131 (1996). 


\bibitem{o.col} Osada,~H., \ttl{Non-collision and collision properties of Dyson's model in infinite dimensions and other stochastic dynamics whose equilibrium states are determinantal random point fields} in Stochastic Analysis on Large Scale Interacting Systems, eds.\ T.\ Funaki and H.\ Osada, Advanced Studies in Pure Mathematics \textbf{39}, 2004, 325-343. 

%\bibitem{o.tp} Osada,~H., \ttl{Tagged particle processes and their non-explosion criteria} J. Math. Soc. Japan, {\bf 62}, No.\ {\bf 3}, 867-894 (2010). 

\bibitem{o.isde} Osada,~H., \ttl{Infinite-dimensional stochastic differential equations related to random matrices} Probability Theory and Related Fields, {\bf 153}, 471-509 (2012). 

\bibitem{o.rm} Osada,~H., \ttl{Interacting Brownian motions in infinite dimensions with logarithmic interaction potentials} Ann. of  Probab. {\bf 41}, 1-49 (2013). 

\bibitem{o.rm2} Osada, H., \ttl{Interacting Brownian motions in infinite dimensions with logarithmic interaction potentials II : Airy random point field} Stochastic Processes and their Applications {\bf 123}, 813-838 (2013). 

\bibitem{o-o.tt} Osada, H., Osada, S., \ttl{Discrete approximations of determinantal point processes on continuous spaces: tree representations and  tail triviality} Journal of Statistical Physics, 170(2018), no.2, 421--435. 

%\bibitem{o-s.08} Osada, H., Shirai, T., \ttl{Variance of the linear statistics of the Ginibre random point field} RIMS K\^oky\^uroku Bessatsu {\bf B6}, 193--200 (2008). 

%\bibitem{o-t.core} Osada,~H., Tanemura,~H., \ttl{Cores of Dirichlet forms related to Random Matrix Theory} (preprint) {\sf http://arxiv.org/abs/1405.4304}.

\bibitem{o-t.sm} Osada,~H., Tanemura,~H., \ttl{Strong Markov property of determinantal processes with extended kernels} Stochastic Processes and their Applications {\bf 126},  186-208, no. 1,  (2016),  doi:10.1016/j.spa.2015.08.003. 

\bibitem{o-t.tail} Osada,~H., Tanemura,~H., \ttl{Infinite-dimensional stochastic differential equations and tail $ \sigma $-fields} arXiv:1412.8674v5 [math.PR]. 

\bibitem{o-t.airy} Osada,~H., Tanemura,~H., \ttl{Infinite-dimensional stochastic differential equations related to Airy random point fields} arXiv:1408.0632.

\bibitem{PR29} Plancherel, M., Rotach, W., \ttl{Sur les valeurs asymptotiques des polynomes d\'Hermite $H_n(x)=(-1)^{n } e^{x^2/2}\frac{d^{n }}{dx^{n }}(e^{-x^2/2})$} Comment. Math. Helv. {\bf 1}, 227-254 (1929). 

%\bibitem{rrv.airy} Ram\'{i}rez, A., Rider B., Vir\'{a}g, B., \ttl{Beta ensembles, stochastic Airy spectrum, and a diffusion} J. Amer. Math. Soc. {\bf 24} (4), 919-944 (2011)

%\bibitem{RY} Revuz D., Yor M., \ttl{Continuous martingales and Brownian motions} (3rd ed) Springer (1999). 

%\bibitem{ruelle.2}  Ruelle,~D., \ttl{Superstable interactions in classical statistical mechanics} Commun. Math. Phys. {\bf 18} (1970) 127--159. 


\bibitem{ST03} Shirai, T., Takahashi, Y., \ttl{ Random point fields associated with certain Fredholm determinants I:\ fermion, Poisson and boson point process} J. Funct. Anal. {\bf 205}, 414-463 (2003). 

\bibitem{sos.00} Soshnikov, A., \ttl{Determinantal random point fields} Russian Math. Surveys {\bf 55}, 923-975 (2000). 

\bibitem{spohn.dyson} Spohn, H., \ttl{Interacting Brownian particles:
a study of Dyson's model} In: Hydrodynamic Behavior and Interacting Particle Systems, G. Papanicolaou (ed), IMA Volumes in Mathematics and its Applications, {\bf 9}, Berlin: Springer-Verlag, pp. 151-179 (1987). 

%\bibitem{Sze} Szeg\"o, G., \ttl{Orthogonal Polynomials} (4th ed) Amer. Math. Sos. Colloquium Publications,  (1981)

%\bibitem{tane.2} Tanemura,~H., \ttl{A system of infinitely many mutually reflecting Brownian balls in $\R^d $}  Probab.\  Theory Relat.\  Fields {\bf  104} (1996) 399-426. 

\bibitem{tsai.14}  Tsai, Li-Cheng,  \ttl{Infinite dimensional stochastic differential equations for Dyson's model} Probab.\  Theory Relat.\  Fields  (published on line) DOI 10.1007/s00440-015-0672-2  (2015).  

\end{thebibliography}
\end{document}